\documentclass[12pt]{article}
\usepackage{latexsym}
\usepackage{amsfonts}
\usepackage{amssymb}
\input epsf.tex
\parskip=\medskipamount                 
\thicklines                     
\newcommand{\bimn}[7]{\bibitem{#1}#2,
{\em #3},
{ #4}\hspace{0.25em}{\bf
#5}\hspace{0.25em}(#6)\hspace{0.25em}{#7}.}

\def\inbar{\vrule height1.5ex width.4pt depth0pt}
\def\IC{\relax\,\hbox{$\inbar\kern-.3em{\rm C}$}}
\def\IN{\relax{\rm I\kern-.18em N}}
\def\IQ{\relax\,\hbox{$\inbar\kern-.3em{\rm Q}$}}
\def\IR{\relax{\rm I\kern-.18em R}}
\def\ZZ{\relax{\sf Z\kern-.4em Z}}
\def\a{\alpha} \def\b{\beta}    
 \def\l{\lambda} 
   
 \def\cB{{\cal B}}  \def\cD{{\cal D}}
   
 \def\cK{{\cal K}} \def\cL{{\cal L}} 
 \def\cO{{\cal O}}  
 
\newtheorem{theorem}{Theorem}[section]

\newtheorem{corollary}[theorem]{Corollary}
\newtheorem{conjecture}[theorem]{Conjecture}
\newtheorem{lemma}[theorem]{Lemma}

\newtheorem{remark}[theorem]{Remark}

\newtheorem{crl}{Corollary of Conjecture~4.1}

\catcode`\@=11

\newif\if@fewtab\@fewtabtrue

\catcode`\@=11

\newif\if@fewtab\@fewtabtrue

{\count255=\time\divide\count255 by 60
\xdef\hourmin{\number\count255}
\multiply\count255 by-60\advance\count255 by\time
\xdef\hourmin{\hourmin:\ifnum\count255<10 0\fi\the\count255}}
\def\ps@draft{\let\@mkboth\@gobbletwo
    \def\@oddhead{}
    \def\@oddfoot
      {\hbox to 7 cm{\footnotesize {\em Draft of \jobname:} \draftdate
       \hfil}\hskip -7cm\hfil\rm\thepage \hfil}
    \def\@evenhead{}\let\@evenfoot\@oddfoot}

\def\ceqno{\global\@fewtabfalse
    \ifcase\@eqcnt \def\@tempa{& & &}\or \def\@tempa{& &}
      \or \def\@tempa{&}
      \or\def\@tempa{}\fi\@tempa
{\rm(\theequation)}}

\def\aeqno#1{\global\@fewtabfalse
    \ifcase\@eqcnt \def\@tempa{& & &}\or \def\@tempa{& &}
      \or \def\@tempa{&}
      \or\def\@tempa{}\fi\@tempa
{\rm(\theequation,#1)}}

\def\label#1{\ifnum\draftcontrol=1
 \global\def\draftnote{$\scriptstyle #1$}\fi
 \@bsphack\if@filesw {\let\thepage\relax
   \def\protect{\noexpand\noexpand\noexpand}%
\xdef\@gtempa{\write\@auxout{\string
      \newlabel{#1}{{\@currentlabel}{\thepage}}}}}\@gtempa
   \if@nobreak \ifvmode\nobreak\fi\fi\fi
  \@esphack}

\def\alabel#1#2{\label{#1}\global\@fewtabfalse
    \ifcase\@eqcnt \def\@tempa{& & &}\or \def\@tempa{& &}
      \or \def\@tempa{&}
      \or\def\@tempa{}\fi\@tempa
{\hbox to 3cm{\phantom{\rm(\theequation,#2)}
\draftnote \hfil}\hskip -3cm {\rm(\theequation,#2)}}}

\def\clabel#1{\label{#1}\global\@fewtabfalse
    \ifcase\@eqcnt \def\@tempa{& & &}\or \def\@tempa{& &}
      \or \def\@tempa{&}
      \or\def\@tempa{}\fi\@tempa
{\hbox to 3cm{\phantom{\rm(\theequation)}
\draftnote \hfil}\hskip -3cm{\rm(\theequation)}}}

\def\eqnarray{\def\draftnote{{}}\global\@fewtabtrue
\stepcounter{equation}\let\@currentlabel=\theequation
\global\@eqnswtrue
\global\@eqcnt\z@\tabskip\@centering\let\\=\@eqncr
$$\halign to \displaywidth\bgroup\@eqnsel\hskip\@centering\@eqcnt\z@
  $\displaystyle\tabskip\z@{##}$&\global\@eqcnt\@ne
  \hskip 1\arraycolsep \hfil$\displaystyle{##}$\hfil
  &\global\@eqcnt\tw@ \hskip 1\arraycolsep
$\displaystyle\tabskip\z@{##}$
\hfil  \tabskip\@centering&\global\@eqcnt\thr@@\llap{##}\tabskip\z@
\cr}

\def\endeqnarray{\@@eqncr\egroup
      \global\advance\c@equation\m@ne$$\global\@ignoretrue}

\def\@eqnnum{\hbox to 3cm{\phantom{\rm(\theequation)} \draftnote
                         \hfil}\hskip -3cm {\rm(\theequation)}}

\def\@@eqncr{\let\@tempa\relax
    \ifcase\@eqcnt \def\@tempa{& & &}\or \def\@tempa{& &}
      \or \def\@tempa{&}
      \or\def\@tempa{}
\fi\@tempa
\if@eqnsw
\if@fewtab\@eqnnum\fi
\stepcounter{equation}\fi\global
\@eqnswtrue\global\@eqcnt\z@\global\@fewtabtrue\cr}

\def\draftcite#1{\ifnum\draftcontrol=1#1\else{}\fi}

\def\@lbibitem[#1]#2{\item{}\hskip -3cm \hbox to 2cm
{\hfil$\scriptstyle\draftcite{#2}$}\hskip
1cm[\@biblabel{#1}]\if@filesw
     {\def\protect##1{\string ##1\space}\immediate
      \write\@auxout{\string\bibcite{#2}{#1}}}\fi\ignorespaces}

\def\@bibitem#1{\item\hskip -3cm \hbox to 2cm
{\hfil $\scriptstyle\draftcite{#1}$}\hskip 1cm
\if@filesw \immediate\write\@auxout
       {\string\bibcite{#1}{\the\value{\@listctr}}}\fi\ignorespaces}

\def\nsection#1{\section{#1}\setcounter{equation}{0}}

\def\draftdate{\number\month/\number\day/\number\year\ \ \ \hourmin }

\global\def\draftcontrol{0}
\catcode`\@=12

\def\theequation{{\thesection.\arabic{equation}}}

\def\qq{\begin{eqnarray}}
\def\qqq{\end{eqnarray}}
\def\rx#1{~(\ref{#1})}
\def\rxw#1{(\ref{#1})}
\def\ex#1{eq.\hspace*{-3pt}\rx{#1}}
\def\eex#1{eqs.\hspace*{-3pt}\rx{#1}}
\def\cx#1{~\cite{#1}}
\def\rw#1{~\ref{#1}}

\newlength{\shiftwidth}
\def\shift#1{&&\hbox to \shiftwidth{\hfill $\displaystyle#1$}}
\newlength{\sshiftwidth}
\def\sshift#1{\lefteqn{\hbox to
\sshiftwidth{\hfill$\displaystyle#1$}}}

\def\ie{{\it i.e.\ }}

\def\cf{{\it cf.\ }}
\def\rhs{{\it r.h.s.\ }}
\def\lhs{{\it l.h.s.\ }}

\def\deg{ \mathop{{\rm deg}}\nolimits }
\def\p{^{\prime}}

\def\pr#1#2{ \noindent{\em Proof of #1~\ref{#2}.} }
\def\proof{ \noindent{\em Proof.} }

\def\qed{ \hfill $\Box$ }

\def\lrbc#1{ \left( #1 \right) }
\def\lrbs#1{ \left[ #1 \right] }

\def\u#1{ \underline{#1} }

\def\um{ {\u{m}} }

\def\sjoN{ \sum_{j=1}^N }

\def\pjoN{ \prod_{j=1}^N }

\def\qa{ q^{\a} }

\def\Ss{ S^* }

\def\Cz{ \IC[\zoN] }

\def\ZZ{ \mathbb{Z} }
\def\IQ{ \mathbb{Q} }
\def\IC{ \mathbb{C} }

\def\qbezier{\bezier{120}}
\def\DottedCircle{
\bezier{4}(0.966,-0.259)(1.04,0)(0.966,0.259)
\bezier{4}(0.966,0.259)(0.897,0.518)(0.707,0.707)
\bezier{4}(0.707,0.707)(0.518,0.897)(0.259,0.966)
\bezier{4}(0.259,0.966)(0,1.04)(-0.259,0.966)
\bezier{4}(-0.259,0.966)(-0.518,0.897)(-0.707,0.707)
\bezier{4}(-0.707,0.707)(-0.897,0.518)(-0.966,0.259)
\bezier{4}(-0.966,0.259)(-1.04,0)(-0.966,-0.259)
\bezier{4}(-0.966,-0.259)(-0.897,-0.518)(-0.707,-0.707)
\bezier{4}(-0.707,-0.707)(-0.518,-0.897)(-0.259,-0.966)
\bezier{4}(-0.259,-0.966)(0,-1.04)(0.259,-0.966)
\bezier{4}(0.259,-0.966)(0.518,-0.897)(0.707,-0.707)
\bezier{4}(0.707,-0.707)(0.897,-0.518)(0.966,-0.259)
}
\def\Endpoint[#1]{
\ifcase#1
\put(1,0){\circle*{0.15}}
\or\put(0.866,0.5){\circle*{0.15}}
\or\put(0.5,0.866){\circle*{0.15}}
\or\put(0,1){\circle*{0.15}}
\or\put(-0.5,0.866){\circle*{0.15}}
\or\put(-0.866,0.5){\circle*{0.15}}
\or\put(-1,0){\circle*{0.15}}
\or\put(-0.866,-0.5){\circle*{0.15}}
\or\put(-0.5,-0.866){\circle*{0.15}}
\or\put(0,-1){\circle*{0.15}}
\or\put(0.5,-0.866){\circle*{0.15}}
\or\put(0.866,-0.5){\circle*{0.15}}
\fi}
\def\Arc[#1]{
\thicklines         
\ifcase#1
\bezier{25}(0.966,-0.259)(1.04,0)(0.966,0.259)
\or
\bezier{25}(0.966,0.259)(0.897,0.518)(0.707,0.707)
\or
\bezier{25}(0.707,0.707)(0.518,0.897)(0.259,0.966)
\or
\bezier{25}(0.259,0.966)(0,1.04)(-0.259,0.966)
\or
\bezier{25}(-0.259,0.966)(-0.518,0.897)(-0.707,0.707)
\or
\bezier{25}(-0.707,0.707)(-0.897,0.518)(-0.966,0.259)
\or
\bezier{25}(-0.966,0.259)(-1.04,0)(-0.966,-0.259)
\or
\bezier{25}(-0.966,-0.259)(-0.897,-0.518)(-0.707,-0.707)
\or
\bezier{25}(-0.707,-0.707)(-0.518,-0.897)(-0.259,-0.966)
\or
\bezier{25}(-0.259,-0.966)(0,-1.04)(0.259,-0.966)
\or
\bezier{25}(0.259,-0.966)(0.518,-0.897)(0.707,-0.707)
\or
\bezier{25}(0.707,-0.707)(0.897,-0.518)(0.966,-0.259)
\fi}
\def\DottedArc[#1]{
\ifcase#1
\bezier{4}(0.966,-0.259)(1.04,0)(0.966,0.259)
\or
\bezier{4}(0.966,0.259)(0.897,0.518)(0.707,0.707)
\or
\bezier{4}(0.707,0.707)(0.518,0.897)(0.259,0.966)
\or
\bezier{4}(0.259,0.966)(0,1.04)(-0.259,0.966)
\or
\bezier{4}(-0.259,0.966)(-0.518,0.897)(-0.707,0.707)
\or
\bezier{4}(-0.707,0.707)(-0.897,0.518)(-0.966,0.259)
\or
\bezier{4}(-0.966,0.259)(-1.04,0)(-0.966,-0.259)
\or
\bezier{4}(-0.966,-0.259)(-0.897,-0.518)(-0.707,-0.707)
\or
\bezier{4}(-0.707,-0.707)(-0.518,-0.897)(-0.259,-0.966)
\or
\bezier{4}(-0.259,-0.966)(0,-1.04)(0.259,-0.966)
\or
\bezier{4}(0.259,-0.966)(0.518,-0.897)(0.707,-0.707)
\or
\bezier{4}(0.707,-0.707)(0.897,-0.518)(0.966,-0.259)
\fi}
\def\Chord[#1,#2]{
\thinlines
\ifnum#1>#2\Chord[#2,#1]
\else\ifnum#1<#2
\ifcase#1
\ifcase#2
\or\qbezier(1,0)(0.516,0.138)(0.866,0.5)
\or\qbezier(1,0)(0.45,0.26)(0.5,0.866)
\or\qbezier(1,0)(0.327,0.327)(0,1)
\or\qbezier(1,0)(0.179,0.311)(-0.5,0.866)
\or\qbezier(1,0)(0.0536,0.2)(-0.866,0.5)
\or\put(1, 0){\line(-2, 0){2}}
\or\qbezier(1,0)(0.0536,-0.2)(-0.866,-0.5)
\or\qbezier(1,0)(0.179,-0.311)(-0.5,-0.866)
\or\qbezier(1,0)(0.327,-0.327)(0,-1)
\or\qbezier(1,0)(0.45,-0.26)(0.5,-0.866)
\or\qbezier(1,0)(0.516,-0.138)(0.866,-0.5)
\fi
\or\ifcase#2\or
\or\qbezier(0.866,0.5)(0.378,0.378)(0.5,0.866)
\or\qbezier(0.866,0.5)(0.26,0.45)(0,1)
\or\qbezier(0.866,0.5)(0.12,0.446)(-0.5,0.866)
\or\qbezier(0.866,0.5)(0,0.359)(-0.866,0.5)
\or\qbezier(0.866,0.5)(-0.0536,0.2)(-1,0)
\or\put(0.866, 0.5){\line(-5, -3){1.73}}
\or\qbezier(0.866,0.5)(0.146,-0.146)(-0.5,-0.866)
\or\qbezier(0.866,0.5)(0.311,-0.179)(0,-1)
\or\qbezier(0.866,0.5)(0.446,-0.12)(0.5,-0.866)
\or\qbezier(0.866,0.5)(0.52,0)(0.866,-0.5)
\fi
\or\ifcase#2\or\or
\or\qbezier(0.5,0.866)(0.138,0.516)(0,1)
\or\qbezier(0.5,0.866)(0,0.52)(-0.5,0.866)
\or\qbezier(0.5,0.866)(-0.12,0.446)(-0.866,0.5)
\or\qbezier(0.5,0.866)(-0.179,0.311)(-1,0)
\or\qbezier(0.5,0.866)(-0.146,0.146)(-0.866,-0.5)
\or\put(0.5, 0.866){\line(-3, -5){1}}
\or\qbezier(0.5,0.866)(0.2,-0.0536)(0,-1)
\or\qbezier(0.5,0.866)(0.359,0)(0.5,-0.866)
\or\qbezier(0.5,0.866)(0.446,0.12)(0.866,-0.5)
\fi
\or\ifcase#2\or\or\or
\or\qbezier(0,1.)(-0.138,0.516)(-0.5,0.866)
\or\qbezier(0,1.)(-0.26,0.45)(-0.866,0.5)
\or\qbezier(0,1.)(-0.327,0.327)(-1,0)
\or\qbezier(0,1.)(-0.311,0.179)(-0.866,-0.5)
\or\qbezier(0,1.)(-0.2,0.0536)(-0.5,-0.866)
\or\put(0, 1){\line(0, -2){2}}
\or\qbezier(0,1.)(0.2,0.0536)(0.5,-0.866)
\or\qbezier(0,1.)(0.311,0.179)(0.866,-0.5)
\fi
\or\ifcase#2\or\or\or\or
\or\qbezier(-0.5,0.866)(-0.378,0.378)(-0.866,0.5)
\or\qbezier(-0.5,0.866)(-0.45,0.26)(-1,0)
\or\qbezier(-0.5,0.866)(-0.446,0.12)(-0.866,-0.5)
\or\qbezier(-0.5,0.866)(-0.359,0)(-0.5,-0.866)
\or\qbezier(-0.5,0.866)(-0.2,-0.0536)(0,-1)
\or\put(-0.5, 0.866){\line(3, -5){1}}
\or\qbezier(-0.5,0.866)(0.146,0.146)(0.866,-0.5)
\fi
\or\ifcase#2\or\or\or\or\or
\or\qbezier(-0.866,0.5)(-0.516,0.138)(-1,0)
\or\qbezier(-0.866,0.5)(-0.52,0)(-0.866,-0.5)
\or\qbezier(-0.866,0.5)(-0.446,-0.12)(-0.5,-0.866)
\or\qbezier(-0.866,0.5)(-0.311,-0.179)(0,-1)
\or\qbezier(-0.866,0.5)(-0.146,-0.146)(0.5,-0.866)
\or\put(-0.866, 0.5){\line(5, -3){1.73}}
\fi
\or\ifcase#2\or\or\or\or\or\or
\or\qbezier(-1,0)(-0.516,-0.138)(-0.866,-0.5)
\or\qbezier(-1,0)(-0.45,-0.26)(-0.5,-0.866)
\or\qbezier(-1,0)(-0.327,-0.327)(0,-1)
\or\qbezier(-1,0)(-0.179,-0.311)(0.5,-0.866)
\or\qbezier(-1,0)(-0.0536,-0.2)(0.866,-0.5)
\fi
\or\ifcase#2\or\or\or\or\or\or\or
\or\qbezier(-0.866,-0.5)(-0.378,-0.378)(-0.5,-0.866)
\or\qbezier(-0.866,-0.5)(-0.26,-0.45)(0,-1)
\or\qbezier(-0.866,-0.5)(-0.12,-0.446)(0.5,-0.866)
\or\qbezier(-0.866,-0.5)(0,-0.359)(0.866,-0.5)
\fi
\or\ifcase#2\or\or\or\or\or\or\or\or
\or\qbezier(-0.5,-0.866)(-0.138,-0.516)(0,-1)
\or\qbezier(-0.5,-0.866)(0,-0.52)(0.5,-0.866)
\or\qbezier(-0.5,-0.866)(0.12,-0.446)(0.866,-0.5)
\fi
\or\ifcase#2\or\or\or\or\or\or\or\or\or
\or\qbezier(0,-1.)(0.138,-0.516)(0.5,-0.866)
\or\qbezier(0,-1.)(0.26,-0.45)(0.866,-0.5)
\fi
\or\ifcase#2\or\or\or\or\or\or\or\or\or\or
\or\qbezier(0.5,-0.866)(0.378,-0.378)(0.866,-0.5)
\fi\fi\fi\fi}
\def\FullChord[#1,#2]{
\Endpoint[#1]
\Endpoint[#2]
\Arc[#1]
\Arc[#2]
\Chord[#1,#2]
}
\def\EndChord[#1,#2]{
\Endpoint[#1]
\Endpoint[#2]
\Chord[#1,#2]
}
\def\Picture#1{
\begin{picture}(2,1)(-1,-0.167)
#1
\end{picture}
}
\def\DottedChordDiagram[#1,#2]{
\Picture{\DottedCircle \FullChord[#1,#2]}
}
\def\u#1{ \underline{#1} }

\def\um{ {\u{m}} }

\def\vx{ \vec{x} }
\def\vy{ \vec{y} }

\def\val{ \vec{\a} }

\def\val{ \vec{\a} }

\def\bx{ \bar{x} }
\def\by{ \bar{y} }

\def\tI{ \tilde{I} }

\def\tb{ \tilde{b} }

\def\pjoN{ \prod_{j=1}^N }

\def\snzi{ \sum_{n=0}^\infty }

\def\snzi{ \sum_{n=0}^\infty }

\def\smzi{ \sum_{m=0}^\infty }

\def\Jbas#1#2#3{ J^{#2}_{#1}(#3;q) }

\def\Pbas#1#2#3{ \Pbasa{#1}{#2;#3} }

\def\Pi#1{ P_{#1} }

\def\cLp{ \cL\p }
\def\cLz{ \cL_{0} }

\def\qa{ q^\a }

\def\hlfv{ {1\over 2} }

\def\val{ \vec{\a} }

\def\yp{ p }

\def\cLp{ \cL\p }

\def\ell{,\ldots,}

\def\hlfv{ {1\over 2} }

\def\snzi{ \sum_{n=0}^\infty }
\def\snoi{ \sum_{n=1}^\infty }

\def\AP{ \Delta_{\rm A} }

\def\APbas#1#2{ \AP(#1;#2) }

\def\Jbas#1#2{ J_{#1}(#2;q) }

\def\empt{ \varnothing }

\def\Pbas#1#2{ P_{#1}(#2) }

\def\pbas#1#2{ p_{#1}(#2) }

\def\Lbas#1#2#3 { L_{#1}(#2;#3) }

\def\cLz{ \cL_0 }

\def\bt{ \mathbf{t} }
\def\bx{ \mathbf{x} }
\def\by{ \mathbf{y} }

\def\cB{ \mathcal{B} }
\def\tcB{ \tilde{\cB} }
\def\tcD{ \tilde{\cD} }

\def\AS{ {\rm AS} }
\def\tcBas{ \tcB_\AS }
\def\tcBs{ \tcB\p }
\def\tcBr{ \tcB\rst }
\def\tcBrz{ \tcB^{({\rm r},0)} }
\def\IHX{ {\rm IHX} }
\def\tcBihx{ \tcB_\IHX }

\def\rst{ ^{({\rm r})} }

\def\tcBasri#1{ \tcBas^{(#1)} }
\def\tcBihxri#1{ \tcBihx^{(#1)} }
\def\tcBasro{ \tcBasri{1} }
\def\tcBihxro{ \tcBihxri{1} }

\def\cchBihxro{ \cchB_\IHX^{(1)} }

\def\tcBasrz{ \tcBasri{0} }
\def\tcBasrii{ \tcBasri{i} }

\def\cchBihxrodi#1{ \cchBihxro(D,#1) }
\def\cchBihxrodm{ \cchBihxrodi{m} }
\def\cchBihxrodo{ \cchBihxrodi{1} }

\def\tcBihxrodi#1{ \tcBihxro(D,#1) }
\def\tcBihxrodm{ \tcBihxrodi{m} }

\def\tcBihxi#1{ \tcBihx^{(#1)} }
\def\tcBihxo{ \tcBihxi{1} }
\def\tcBihxz{ \tcBihxi{0} }
\def\tcBihxii{ \tcBihxi{i} }

\def\cBs{ \cB^* }

\def\sjoN{ \sum_{j=1}^N }

\def\pjoN{ \prod_{j=1}^N }

\def\hfAS{ \hat{f}_{\AS} }
\def\hf{ \hat{f} }
\def\hfs{ \hf^* }

\def\tV{ \tilde{V} }
\def\tVAS{ \tV_{\AS} }
\def\VAS{ V_{\AS} }

\def\tcDAS{ \tcD_{\AS} }
\def\tcDIHX{ \tcD_{\IHX} }

\def\AP{ \Delta_{\rm A} }
\def\APbas#1#2{ \AP(#1;#2) }
\def\APtbas#1{ \APbas{#1}{t} }
\def\APLt{ \APtbas{\cL} }
\def\APKt{ \APtbas{\cK} }
\def\APKqa{ \APbas{\cK}{\qa} }
\def\APKqai#1{ \AP^{#1}(\cK;\qa) }
\def\APKqan{ \APKqai{2n} }

\def\qa{ q^{\a} }

\def\Jbas#1#2{ J_{#1}(#2;q) }
\def\Jt#1{ \Jbas{2}{#1} }
\def\JtL{ \Jt{\cL} }
\def\JtK{ \Jt{\cK} }
\def\JvaK{ \Jbas{\val}{\cK} }

\def\JaK{ \Jbas{\a}{\cK} }

\def\Pbas#1#2#3{ P_{#1}(#2;#3) }
\def\PnKt{ \Pbas{n}{\cK}{t} }
\def\PnKqa{ \Pbas{n}{\cK}{\qa} }

\def\pbas#1#2#3{ p_{#1}(#2;#3) }
\def\pnvaK{ \pbas{n}{\cK}{\val} }

\def\ZZti{ \ZZ[t^{\pm 1}] }
\def\ZZqhi{ \ZZ[q^{\pm1/2}] }

\def\snzi{ \sum_{n=0}^\infty }
\def\smnog{ \sum_{m,n\geq 0\atop m+n\geq 1} }
\def\smnoge{ \sum_{{m\geq 0\atop n\geq m}\atop m+n\geq 1} }
\def\sD{ \sum_{D\in\bD} }
\def\sDc{ \sum_{D\in\bDc} }
\def\sDl{ \sum_{D\in\bD,\,\chi(D)\geq 1} }
\def\sDlo{ \sum_{D\in\bD,\,\chi(D)= 1} }
\def\sDln{ \sum_{D\in\bD,\,\chi(D) = n} }

\def\bDc{ \bD_{\rm c} }

\def\bopnzi{ \bigoplus_{n=0}^\infty }
\def\bopmzi{ \bigoplus_{m=0}^\infty }
\def\bopmnzi{ \bigoplus_{m,n=0}^\infty }

\def\cLp{ \cL_+ }
\def\cLm{ \cL_- }
\def\cLz{ \cL_0 }

\def\unknot{ {\rm unknot} }

\def\span#1{\mathop{{\rm span}}\nolimits(#1)}
\def\Inv#1#2{\mathop{{\rm Inv}_{#1}}\nolimits(#2)}
\def\Inv#1#2{(#2)_{#1}}
\def\Invgd#1{ \Inv{G_D}{#1} }
\def\Invg#1{ \Inv{G}{#1} }

\def\Invgdo#1{ \Inv{\SgroD}{#1} }
\def\Invgdx#1{ \Inv{\SgrD}{#1} }
\def\Invgdix#1#2{ \Inv{\Sgri{#1}}{#2} }

\def\dihx{ \partial_\IHX }
\def\hdihx{ \hat{\partial}_\IHX }

\def\Qalghzi#1{ \Inv{\Sgri{#1}}{\IQ[\Hoz{#1}]} }
\def\Qalghz{ \Qalghzi{D} }
\def\Qalghzo{ \Qalghzi{D_1}}

\def\cord{ \star }

\def\Ho#1{ H^1(#1) }
\def\Hor#1{ \Ho{#1,\IQ} }
\def\Hoz#1{ \Ho{#1,\ZZ} }
\def\Hsal{ {\cal H} }
\def\Hs#1{ \Hsal(#1) }
\def\Hsi#1#2{ \Hsal_{#1}(#2) }
\def\Hsm#1{ \Hsi{m}{#1} }
\def\SHor#1{ S^*\Hor{#1} }
\def\SiHor#1#2{ S^{#1}\Hor{#2} }
\def\SsHor#1{ \SiHor{*}{#1} }
\def\SmHor#1{ \SiHor{m}{#1} }

\def\Hsalo{ {\cal \Hsal^{\cord} } }
\def\Hso#1{ \Hsalo(#1) }
\def\Hsio#1#2{ \Hsal^{\cord}_{#1}(#2) }
\def\Hsmo#1{ \Hsio{m}{#1} }

\def\Hsalx{ {\cal \Hsal^{\xcord} } }
\def\Hsx#1{ \Hsalx(#1) }
\def\Hsix#1#2{ \Hsal^{\xcord}_{#1}(#2) }
\def\Hsmx#1{ \Hsix{m}{#1} }

\def\xcord{}

\def\Pgio#1{ P^{\cord}_{#1} }
\def\PgoD{ \Pgio{D} }

\def\Pgix#1{ P^{\xcord}_{#1} }
\def\PgxD{ \Pgix{D} }

\def\tcDo{ \tcD^{\cord} }
\def\tcDio#1{ \tcD^{\cord}_{#1} }
\def\tcDix#1{ \tcD^{\xcord}_{#1} }

\def\tcDASo{ \tcDio{\AS} }
\def\tcDIHXo{ \tcDio{\IHX} }
\def\tcDIHXx{ \tcDix{\IHX} }

\def\tcDbd{ \tcD_{\bD} }

\def\dsrdr#1#2{ \#(#1,#2) }
\def\dsd{ \dsrdr{D_1}{D_2} }

\def\Sgri#1{ G_{#1} }
\def\SgrD{ \Sgri{D} }
\def\Sgro{ G^{\cord} }
\def\Sgroi#1{ \Sgro_{#1} }
\def\SgroD{ \Sgroi{D} }

\def\lini#1{ \chi(#1)+1 }
\def\linD{ \lini{D} }
\def\fbs{ f_1,\cdots,f_{\linD} }
\def\efbs{ e^{f_1},\cdots,e^{f_{\linD}} }
\def\tbs{ t_1,\cdots,t_{\linD} }
\def\xbs{ x_1,\cdots,x_{\linD} }

\def\xot{ x_1,x_2 }
\def\exot{ e^{x_1},e^{x_2} }

\def\tot{ t_1,t_2 }

\def\zth{ I_\theta }
\def\zthK#1{ \zth(\cK;#1) }
\def\zthKx{ \zthK{\xot} }

\def\zths{ I_\theta^* }
\def\zthsK#1{ \zths(\cK;#1) }
\def\zthsKt{ \zthsK{\tot} }

\def\pth{ p_\theta }
\def\pthK#1{ \pth(\cK;#1) }
\def\pthKt{ \pthK{\tot} }
\def\pthKex{ \pthK{\exot} }
\def\pthKef{ \pthK{\fDeot} }

\def\APK#1{ \APbas{\cK}{#1} }
\def\APKe#1{ \APK{e^{#1}} }
\def\APKex#1{ \APKe{x_{#1}} }

\def\APKti#1{ \APK{t_{#1}} }

\def\fD#1#2{ f_{#1,D_{#2}} }
\def\fDot#1{ \fD{1}{#1},\fD{2}{#1} }
\def\fDeot{ e^{\fD{1}{1}},e^{\fD{2}{1}} }
\def\efD#1{ e^{f_{#1,D_1}} }

\def\Hsri#1#2{ H^{({\rm exp})}(#1,#2) }
\def\HsQ#1{ \Hsri{#1}{\IQ} }

\def\oHo#1{ H_1(#1) }
\def\oHor#1{ \oHo{#1,\IQ} }

\def\rnk#1{ |#1| }

\def\dego{ \deg_1 }
\def\degt{ \deg_2 }

\def\hxA{ \hat{\xA} }
\def\xA{ A }

\def\sass{ \sum_{c\in\bfS} }

\def\Co{ C_1 }

\def\Cos{ C_1^* }
\def\SmCos{ S^m\Cos }

\def\chCos{ \chC_1^* }
\def\chC{ \check{C} }

\def\Cz{ C_0 }

\def\xdel{ \partial }

\def\tcBmd{ \tcB_m(D) }
\def\tcBmdz{ \tcB_m(\Dz) }
\def\cBmd{ \cB_m(D) }

\def\cchB{ \check{\cB} }
\def\cchBmdp{ \cchB_m(D) }
\def\cchBmdpz{ \cchB_m(\Dz) }
\def\cchBodp{ \cchB_1(D) }

\def\cBmd{ \cB_m(D) }

\def\cDs{ \cD^* }

\def\bfS{ \mathbf{S} }
\def\tbfS{ \tilde{\bfS} }

\def\Prbas#1{ P_{#1} }
\def\Prl{ \Prbas{\l} }
\def\Ph{ \Prbas{\gh} }

\def\Expalg#1{ \mathop{{\rm Exp}}(#1) }
\def\Exp{ {\rm Exp} }

\def\SmHord{ S^m\Hor{D} }

\def\Ps{ P_{\rm S} }

\def\fg#1{Fig.~\ref{#1}}
\def\tb#1{Table~\ref{#1}}

\def\Nz{ N_0 }

\def\IB{ I^{\cB} }
\def\ID{ I^{\cD} }
\def\IO{ I^{\Omega} }

\def\Il{ I^{({\rm log})} }
\def\tIl{ \tI^{({\rm log})} }
\def\tI{ \tilde{I} }

\def\IlKi#1#2{ \Il(\cK,#1;#2) }

\def\IQtot{ \IQ[t_1^{\pm 1},t_2^{\pm 1} ] }
\def\ZZtot{ \ZZ[t_1^{\pm 1},t_2^{\pm 1} ] }
\def\ZZtok{ \ZZ[t_1^{\pm 1},\ldots,t_{\rg}^{\pm 1}] }
\def\IQtorg{ \IQ[t_1^{\pm 1},\ldots,t_{\rg}^{\pm 1}] }
\def\IQtbs{ \IQ[t_1^{\pm 1},\ldots,t_{\chi(D)+1}^{\pm 1}] }

\def\Jl#1{ J^{\rm (log) }_{#1} }
\def\JlKi#1{ \Jl{#1}(\cK;\val) }
\def\JlKn{ \JlKi{n} }
\def\JlKo{ \JlKi{1} }

\def\plgi#1{ p^{\rm (log)}_{#1} }
\def\plgn{ \plgi{n} }

\def\bD{ \mathbf{D} }

\def\hb{ \hbar }

\def\hOm{ \hat{\Omega} }

\def\gg{ \mathfrak{g} }
\def\gh{ \mathfrak{h} }

\def\val{ \vec{\a} }
\def\vr{ \vec{\rho} }
\def\vx{ \vec{x} }
\def\vy{ \vec{y} }

\def\fabco{ f_{ab}^{\;\;\;c} }

\def\Dz{ D_0 }

\def\pDval{ \yp(D,\val) }
\def\pD{ \yp(D) }
\def\pval{ \yps(D,\val) }
\def\pxval{ \yp(x,\val) }
\def\pcDval{ \yp_c(D,\val) }
\def\pcxval{ \yp_c(x,\val) }
\def\pcDzval{ \yp_c(\Dz,\val) }
\def\pDxval{ \ypD(\hxA\,x,\val) }
\def\pDyval{ \ypD(y,\val) }

\def\pcDz{ \ypw_c(\Dz) }
\def\pcD{ \ypw_c(D) }
\def\pcDo{ \ypw_c(D_1) }
\def\pccir{ \ypw_c(\crcl) }

\def\pjov#1{ \prod_{j=1}^{#1} }

\def\pjok{ \prod_{j=1}^k }

\def\pjoth{ \prod_{j=1}^3 }

\def\Dg{ \Delta_\gg }
\def\tok{ t_1,\ldots,t_{\rg} }
\def\qaok{ q^{\val\cdot\l_1},\ldots,q^{\val\cdot\l_{\rg}} }
\def\qaot{ q^{\val\cdot\l_1},q^{\val\cdot\l_2} }

\def\qaic#1{ q^{\val\cdot\lci{#1}} }
\def\qaotc{ \qaic{1},\qaic{2} }
\def\qalot{ q^{\val\cdot(\l_1 + \l_2)} }

\def\aockh{ \hb(\val\cdot\lci{1}),\ldots, \hb(\val\cdot\lci{\rg}) }

\def\pjochi#1{ \prod_{j=1}^{3\chi(#1)} }
\def\pjochD{ \pjochi{D} }

\def\dqval{ d_{q,\val} }
\def\dval{ d_{\val} }
\def\Cqg{ C_{q,\gg} }

\def\ypw{ w }

\def\kg{ k }
\def\rg{ r }

\def\lgdval{ \log\lrbc{ \JvaK\over \Jvau} }

\def\lgdval{ \log\lrbc{ \JvaK / d_{\val}} }

\def\FoK#1{ F_1(\cK;#1) }

\def\yp{ w^{\Omega} }
\def\yps{ w }
\def\ypD{ w^{\cD} }

\def\dval{ d_{\val} }

\def\lcj{ \l_{c(j)} }
\def\lci#1{ \l_{c(#1)} }

\def\ecval{ e_{c,\val} }

\def\crcl{ {\rm circle} }

\begin{document}

\begin{titlepage}
\vfill
\begin{center}

{\large \bf A rationality conjecture about Kontsevich integral of knots and its
implications to the structure of the colored Jones polynomial}\\

\bigskip

\bigskip
\centerline{L. Rozansky\footnote{
This work was supported by NSF Grant DMS-9704893}
}



\centerline{\em Department of Mathematics, University of North Carolina}
\centerline{\em CB \#3250, Phillips Hall}
\centerline{\em Chapel Hill, NC 27599}
\centerline{{\em E-mail address:} {\tt rozansky@math.unc.edu}}

\vfill
{\bf Abstract}

\end{center}
\begin{quotation}

We formulate a conjecture about the structure of Kontsevich integral
of a knot. We describe its value in terms of the generating
functions for the numbers of external edges attached to closed
3-valent diagrams. We conjecture that these functions are rational
functions of the exponentials of their arguments, their denominators
being the powers of the Alexander-Conway polynomial. This conjecture implies
the existence of an expansion of a colored Jones (HOMFLY) polynomial in powers
of $q-1$ whose coefficients are rational functions of $q^\a$ ($\a$ being the
color assigned to the knot). We show how to derive the first Kontsevich
integral polynomial associated to the $\theta$-graph from the rational
expansion of the colored $SU(3)$ Jones polynomial.


\end{quotation} \vfill \end{titlepage}

\pagebreak

\nsection{Introduction}
\label{s1}
\hyphenation{Re-she-ti-khin}
\hyphenation{Tu-ra-ev}

The quantum invariants of knots, links and 3-manifolds, such as the
Jones polynominal and the Witten-Reshetikhin-Turaev invariant, were
discovered about 10 years ago. However, their interpretation in terms
of classical 3-dimensional topology still remains a mystery.

Let us compare the skein relation definition of the Jones polynomial
to that of a much older Alexander-Conway polynomial. The
single-variable Alexander-Conway polynomial $\APLt\in\ZZti$ is a
unique invariant of links in $S^3$ which satisfies the following two
properties. First, the normalization condition:
\qq
\APtbas{\unknot} = 1.
\label{1.1}
\qqq
Second, if $\cLp$, $\cLm$ and $\cLz$ are three links whose regular
projection on a plane is the same except at one spot (see \fg{f1}),
\begin{figure}[hbt]
\leavevmode \centerline{\epsfbox{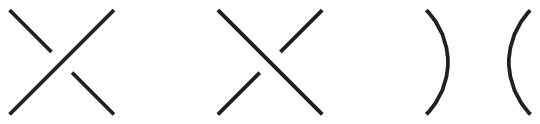}} \caption{The links
$\cLp$, $\cLm$ and $\cLz$}
\label{f1}
\end{figure}
then
\qq
\APtbas{\cLp} - \APtbas{\cLm} = (t^{1/2}-t^{-1/2})\,\APtbas{\cLz}.
\label{1.2}
\qqq
\def\pio{ \pi_1 }
This definition is purely combinatorial and it is a bit unnatural
from the 3-dimensional point of view, since it requires a projection
of a link. However, there exist alternative definitions of the
Alexander-Conway polynomial of a knot $\cK$ which are purely topological.
One derives $\APKt$ from the structure of the knot group
$\pi_1(S^3\setminus \cK)$, and the variable $t$ represents the action of the
homology $\pio/[\pio,\pio]$ onto the quotient
$[\pio,\pio]/[[\pio,\pio],[\pio,\pio]]$, where $\pio$ is the group of the knot
($\pio = \pio(S^3\setminus \cK)$). The other definition
relates the Alexander polynomial to the
Reidemeister torsion of a local system in the knot
complement, the variable $t$ being the twist acquired by that system
along the meridian of $\cK$. From both definitions of $\APKt$ it is
clear that $t$ is intimately related to the meridian of $\cK$.

The Jones polynomial of links $\JtL\in\ZZqhi$ can also be defined by
skein relations. It is the unique invariant which satisfies the
following two properties: the normalization condition
\qq
\Jt{\unknot} = q^{1/2} + q^{-1/2}
\label{1.3}
\qqq
and the skein relation
\qq
q\,\Jt{\cLp} - q^{-1}\,\Jt{\cLm} = (q^{1/2} - q^{-1/2})\,\Jt{\cLz},
\label{1.4}
\qqq
where the links $\cLp$, $\cLm$ and $\cLz$ are the same as those in
\ex{1.2}. Despite an obvious similarity between \eex{1.1},\rx{1.2}
and \eex{1.3},\rx{1.4}, there exists no intrepretation of $\JtL$ in
terms of the ``classical'' objects of 3-dimensional topology, such as
the fundamental group of the knot complement. In particular, there is
no indication that the variable $q$ has any connection to the
meridian of $\cK$.

A new hope for a topological interpretation of $\JtL$ emerged when
J.~Birman, X.-S.~Lin and D.~Bar-Natan discovered that both the
Alexander-Conway and Jones polynomials are packed with Vassiliev
invariants. Consider the expansions
\qq
&
\APKt = \snzi a_n(\cK)\,(t-1)^n,
\label{1.5}\\
&
\JtK = \snzi b_n(\cK)\,(q-1)^n.
\label{1.6}
\qqq
It is not hard to see from the skein relations\rx{1.2} and\rx{1.4}
that the coefficients $\a_n(\cK)$ and $\b_n(\cK)$ are Vassiliev
invariants of degree $n$. However, Vassiliev invariants by definition
are related to the topology of ``the space of all maps
$S^1\longrightarrow S^3$, rather than to the topology of knots
themselves. The latter relation is still missing, although some bits
of it are known, such as the relation between the tree Vassiliev
invariants and Milnor's linking numbers (see \cx{HM} and references therein).

By looking at \ex{1.5} we may say that the Alexander-Conway
polynomial presents a way of assembling some Vassiliev invariants of
knots into a polynomial which has a clear interpretation in terms of
the classical 3-dimensional topology. At the same time, the Jones
polynomial assembles some other Vassiliev invariants into another
polynomial whose topological origin is rather obscure. Therefore one
may wonder if there is a way of reassembling all Vassiliev invariants
into the polynomials which would be similar to the Alexander-Conway
polynomial rather than to the Jones polynomial in terms of their
topological interpretation.

In Sections\rw{gs} and\rw{s3} we present
an algorithm of assembling Vassiliev invariants coming from the Kontsevich
integral of
a knot into a sequence of functions of a variable $t$.
In Section\rw{s4}
we conjecture
that these functions are rational: their denominators are powers of
the Alexander-Conway polynomial of $t$ while their numerators are new
polynomial invariants of knots. Since these new polynomials depend on
the same variable $t$, we expect them to have a topological
interpretation in which, similarly to the case of the
Alexander-Conway polynomial, $t$ will also be related to the meridian
of a knot.

Since the first version of this paper was written and reported, Andrew Kricker
has proved the rationality conjecture in his paper\cx{Kr}.

Kontsevich integral is related to the colored Jones (HOMFLY) polynomial of the
knot through the application of a Lie algebra weight system. In Section\rw{s5}
we explain how to apply this weight system to the `repackaged' Vassiliev
invariants. Then we show how the rational structure of Kontsevich integral
appears as a rational the Jones polynomial. In Section\rw{s6} we use these
results to extract the first non-trivial knot polynomial related to the
$\theta$-graph from the expansion of the $SU(3)$ colored Jones polynomial. In
the Appendix we present a table of these `2-loop' polynomials for knots with up
to 7 crossings.

\nsection{Graph spaces}
\label{gs}

We are going to define an algebra $\cD$ based on 3-valent graphs, but first let
us recall the definition of the algebra $\cB$ of (1,3)-valent graphs related
to Vassiliev invariants of a knot.
Each 3-valent vertex of a graph
is endowed with a cyclic ordering of 3 egdes attached to it. When we draw a
picture of a graph, we assume that this ordering is counterclockwise.

A graph $D$ has 2 degrees. They are defined as
\qq
\dego (D) & = & \#\mbox{1-vertices},
\label{1.7*}
\\
\degt (D) & = & \#\mbox{chords} - \#\mbox{3-vertices} =
\chi(D) + \dego (D),
\label{1.7}
\qqq
where $\chi(D)$ is the Euler characteristic of $D$
(more precisely, $\chi(D)$ denotes the Euler characteristic with the
\emph{opposite} sign).

Let $\tcB_{m,n}$ be a
formal vector space (over $\IC$) whose basis elements are in a
one-to-one correspondence with (1,3)-valent graphs of degrees $m$ and $n$
\qq
\tcB_{m,n} = \span{D\,|\,\dego(D)=m,\,\degt(D)=n}.
\qqq
Together all such spaces form a graded space $\tcB$
\qq
\tcB = \bopnzi \tcB_n, \qquad\mbox{where}\;\;
\tcB_n =
\bopmzi\tcB_{m,n}.
\label{1.8}
\qqq
The space $\tcB$ has two important subspaces:
$\tcBas$ and $\tcBihx$.
$\tcBas$ is spanned by the sums $D_1+D_2$, where $D_1$ and $D_2$ are the same
graphs except that they have different cyclic orders
at one 3-valent vertex:
%
\qq
\tcBas = \span{D_1+D_2\quad\mbox{for all pairs $D_1,D_2$}}.
\label{1.9}
\qqq
In order to define $\tcBihx$, consider a space $\tcBs$ whose basis
vectors are graphs with 1-valent and 3-valent vertices and exactly
one 4-valent vertex. We define a linear map
$\dihx:\,\tcBs\longrightarrow \tcB$ by its action on the individual
graphs of $D\in\tcBs$
\qq
\dihx:\; D\mapsto D_1-D_2+D_3,
\label{1.10}
\qqq
where all four graphs $D,D_1,D_2,D_3$ are the same except at one
spot, where they differ according to \fg{f2}.
\begin{figure}[hbt]
\centerline{\epsfbox{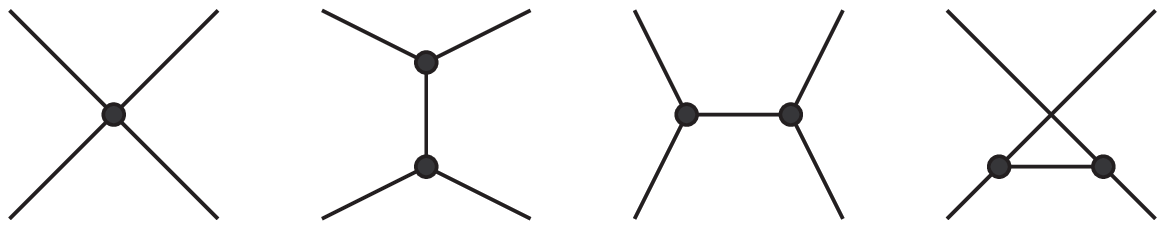}} \caption{The graphs $D$, $D_1$,
$D_2$ and $D_3$}
\label{f2}
\end{figure}
Then we define the second subspace $\tcBihx\subset\tcB$ as the
image of $\dihx$.

Now we introduce a space
\qq
\cB = \tcB / (\tcBas + \tcBihx).
\label{1.11}
\qqq
%
Since the graphs $D_1,D_2$ of\rx{1.9} and the graphs $D_1,D_2,D_3$
of\rx{1.10} have the
same degrees\rx{1.7*},\rxw{1.7} among themselves,
then both subspaces $\tcBas$ and $\tcBihx$ respect the
gradings\rx{1.8} and as a result the space $\cB$ is also graded
\qq
\cB = \bopnzi \cB_n, \qquad\mbox{where}\;\;
\cB_n =
\bopmzi\cB_{m,n}.
\label{1.12}
\qqq
It is well-known that the dual space $\cBs$ is isomorphic to the
space of all Vassiliev invariants of knots, and the grading
$\cBs=\bopnzi \cBs_n$ corresponds to the grading of Vassiliev
invariants.

The space $\cB$ can be endowed with a commutative
algebra structure. We define the product
of two graphs in $\tcB$ as their disjoint union. It is easy to see that this
product respects the gradings\rx{1.7*},\rx{1.7} and that the subspace
$\tcBas + \tcBihx$ is the ideal in algebra $\tcB$. Therefore, the quotient
space $\cB$ is also an algebra.

We are going to introduce another algebra $\cD$ which is isomorphic to
$\cB$. This construction has been known to some people\cx{pr}. It
appeared as an attempt to better understand the structure of $\cB$
and, in particular, to evaluate the dimension of the spaces
$\cB_n$. I am especially indebted to A.~Vaintrob for illuminating
discussions on the structure of $\cD$. I introduce
the algebra $\cD$ in order to formulate a conjecture about the structure
of Kontsevich integral, which was motivated by the study
of the Melvin-Morton expansion of the colored Jones polynomial as it
comes of $R$-matrix expression and which has now been proved by
A.~Kricker\cx{Kr}.

We begin by defining a bigger space $\tcD$.
Let $D$ be a graph with
3-valent vertices and no 1-valent vertices. We think of this graph as
a $CW$-complex and consider a space of its rational cohomologies
$\Hor{D}$. Let $\SgrD$ be the group of symmetry of a graph $D$ (it maps
3-vertices to 3-vertices and edges to edges) and let $\SgroD\subset\SgrD$ be
its subgroup which preserves the cyclic order of the edges at the vertices.
$\SgrD$ acts naturally on $\Hor{D}$ and this group action can be extended to
the symmetric algebra $\SHor{D}$. We denote by $\Hso{D}$ the $\SgroD$-invariant
part of the latter space:
%
\qq
\Hso{D} = \bopmzi \Hsmo{D},\qquad
\Hsio{m}{D} = \Invgdo{\SmHor{D}},
\label{1.13}
\qqq
while $\PgoD$ is the corresponding projector
\qq
\PgoD:\,\SHor{D} \longrightarrow \Hso{D},\qquad
\PgoD(x) = {1\over \rnk{\SgroD}} \sum_{g\in \SgroD} g(x),
\label{1.14}
\qqq
where $\rnk{\SgroD}$ denotes the number of elements in $\SgroD$.
Now we define a linear space $\tcDo$ as
\qq
\tcDo = \bopmnzi \tcDio{m,n},\qquad\mbox{where}\;\;
\tcDio{m,n} = \bigoplus_{D:\,\chi(D)=n} \Hsmo{D}.
\label{1.15}
\qqq

The space $\tcDo$ has an associative, commutative
algebra structure. First, note that for a disjoint
union $D_1\cup D_2$ of two graphs $D_1,D_2$
\qq
\Hor{D_1\cup D_2} = \Hor{D_1} \oplus \Hor{D_2}
\qqq
and therefore
\qq
\SHor{D_1\cup D_2} = \SHor{D_1} \otimes \SHor{D_2}
\qqq
as algebras. The latter equation allows us to define a product of two elements
$x_i\in\Hor{D_i}$, $i=1,2$ as a projection of their tensor product in
$\SHor{D_1\cup D_2}$
\qq
x_1 x_2 = \Pgio{D_1\cup D_2} (x_1\otimes x_2) \in \Hso{D_1\cup D_2}.
\label{1.15*}
\qqq
If the graphs $D_1$, $D_2$ do not have isomorphic connected components, then
$\Sgroi{D_1\cup D_2} = \Sgroi{D_1}\times \Sgroi{D_2}$
and the projector in \ex{1.15*} may
be omitted: $x_1 x_2 = x_1\otimes x_2$. The commutativity of the
product\rx{1.15*} is obvious. Associativity follows from a relation
\qq
(x_1 x_2) x_3 = x_1 (x_2 x_3) = \Pgio{D_1\cup D_2 \cup D_3}
(x_1\otimes x_2\otimes x_3)
\qqq
Finally, since the product\rx{1.15*}
respects both gradings\rx{1.15}, then $\tcDo$ is a graded algebra.

Next, we define the subspace $\tcDASo\subset\tcDo$
which comes from the change of
cyclic order at 3-valent vertices.
The definition of the symmetric algebra $S^*\Hor{D}$ is independent of this
cyclic order. Therefore if we take a graph $D_1$ and
change the cyclic order at one of its vertices, thus producing a new graph
$D_2$, then there is a natural isomorphism of cohomologies
\qq
\hfAS:\; \Hor{D_1}\longrightarrow \Hor{D_2},
\label{1.16}
\qqq
because $D_2$ was constructed in such a way that there is a natural isomorphism
between $D_1$ and $D_2$ as $CW$-complexes (generally, there could be more than
one isomorphism due to the symmetry group $G_{D_1}$). The isomorphism\rx{1.16}
can be extended to an isomorphism of symmetric algebras
\qq
\hfAS:\; \SHor{D_1}\longrightarrow \SHor{D_2},
\label{1.17}
\qqq
%
Let $\tVAS$ be the graph of this map
\qq
\tVAS = \{ (x,y) |\; y = \hfAS(x) \}\subset \SHor{D_1}\oplus
\SHor{D_2}.
\label{1.18}
\qqq
We denote by $\VAS$ its projection onto $\Hso{D_1}\oplus \Hso{D_2}$
\qq
\VAS = \Pgio{D_1}\Pgio{D_2}(\tVAS) \subset \Hso{D_1}\oplus \Hso{D_2}.
\label{1.19}
\qqq
We define the subspace $\tcDASo\subset\tcDo$ as the sum of all the spaces $\VAS$
for all 3-valent diagrams $D_1$ and all choices of vertices of $D_1$ where we
change the orientation. It is easy to check that $\tcDASo$ is an ideal: for any
element $x\in \SHor{D_1}$, and for any element $y\in\Hso{D_3}$
\qq
(\Pgio{D_1}(x) + \Pgio{D_2} \hfAS(x))\, y & = &
\Pgio{D_1\cup D_3}(\Pgio{D_1}(x)\otimes y) +
\Pgio{D_2\cup D_3}(\Pgio{D_2}\hfAS(x)\otimes y)
\nonumber\\
& = &
\Pgio{D_1\cup D_3}(x\otimes y) + \Pgio{D_2\cup D_3} \hfAS(x\otimes y),
\qqq
because obviously
\qq
\Pgio{D_i\cup D_3} (\Pgio{D_i}\otimes I) = \Pgio{D_i\cup D_3},\qquad i=1,2.
\qqq
and
\qq
\hfAS(x\otimes y) = \hfAS(x) \otimes y,
\qqq
where in the \lhs $\hfAS$ comes from the change of cyclic order at a vertex in
the whole graph $D_1\cup D_3$.

Finally, we define a subspace $\tcDIHXo\subset\tcDo$. Let $D$ be a graph with
3-valent vertices and exactly one 4-valent vertex, and with fixed cyclic
order at every vertex.
By adding an extra edge to $D$, we ``resolve'' the 4-valent vertex in
3 different ways, thus converting $D$ into one of the 3-valent graphs
$D_1,D_2,D_3$ of \fg{f2}. A removal of this extra edge generates 3
natural maps of rational homologies
\qq
\hf_i:\; \oHor{D_i} \longrightarrow \oHor{D},\qquad i=1,2,3.
\label{1.20}
\qqq
We extend the dual maps
$\hfs_i:\; \Hor{D}\longrightarrow\Hor{D_i}$ as
algebra homomorphisms
\qq
\hfs_i:\; \SHor{D}\longrightarrow\SHor{D_i},\qquad i=1,2,3.
\label{1.21}
\qqq
We define the map
$\hdihx:\; \SHor{D} \longrightarrow \bigoplus_{i=1}^3\Hso{D_i}$
by the formula (\cf \ex{1.10})
\qq
\hdihx = \sum_{i=1}^3 (-1)^{1+i} \Pgio{D_i} \hfs_i.
\label{1.22}
\qqq
The subspace $\tcDIHXo$ is the sum of the images of all the maps
$\hdihx$ for all the graphs $D$. It is easy to check that similarly to
$\tcDASo$, $\tcDIHXo$ is also an ideal in $\tcDo$.

Now we define the quotient space
\qq
\cD = \tcDo / (\tcDASo + \tcDIHXo).
\label{1.23}
\qqq
Since the graphs $D_1$, $D_2$ of \rxw{1.16} and $D$, $D_1$, $D_2$, $D_3$
of \rxw{1.20} have the same Euler characteristic among themselves and since the
maps\rx{1.17} and\rx{1.21} preserve the grading of symmetric algebras, then
$\cD$ is a graded algebra:
%
\qq
\cD = \bopmnzi \cD_{m,n},
\label{1.24}
\qqq
where the spaces $\cD_{m,n}$ are the quotients of
the spaces $\tcD_{m,n}$.
\qq
\cD_{m,n} = \tcD_{m,n} / \cD_{m,n}\cap (\tcDAS + \tcDIHX).
\label{1.24*}
\qqq
%

This description of the algebra $\cD$
makes it easy to establish its
isomorphism with the algebra $\cB$, but there exists a slightly different
description of $\cD$ which suits better for the formulation of our conjecture
about the structure of Kontsevich integral. Recall that $\SgrD$ denotes the
full symmetry group of a 3-valent graph $D$ (including the maps which do not
preserve the cyclic order at the vertices). As we have mentioned, $\SgrD$ acts
naturally on $\SHor{D}$. We modify this action by multiplying the action of an
element $g\in\SHor{D}$ by $(-1)^{|g|}$, where $|g|$ denotes the number of
vertices of $D$ whose cyclic order is changed by $g$. Now instead of\rx{1.13}
we define
\qq
\Hsx{D} = \bopmzi \Hsmx{D},\qquad
\Hsix{m}{D} = \Invgdx{\SmHor{D}},
\label{1.13*1}
\qqq
while $\PgxD$ is the corresponding projector
\qq
\PgxD:\,\SHor{D} \longrightarrow \Hsx{D},\qquad
\PgxD(x) = {1\over \rnk{\SgrD}} \sum_{g\in \SgrD} g(x),
\label{1.14*1}
\qqq
Let $\bD$ be a set of all 3-valent graphs with a particular cyclic order of
edges at vertices chosen for every graph (so that each isomorphism class of
3-valent graphs is represented in $\bD$ exactly once). Define
\qq
\tcD = \bopmnzi \tcDix{m,n},\qquad\mbox{where}\;\;
\tcDix{m,n} = \bigoplus_{D:\,\chi(D)=n} \Hsmx{D}
\label{1.14*2}
\qqq
(\cf \ex{1.15}).
%
%
If we choose a different set $\bD\p$, then there is a natural isomorphism
between $\tcDbd$ and $\tcD_{\bD\p}$. Namely, if $D_1\in \bD$ and
$D_2\in \bD\p$ represent the same 3-valent graph (but possibly with different
cyclic orders), then we identify the spaces $\Hsx{D_1}$ and $\Hsx{D_2}$ by an
identity map with an extra sign factor $(-1)^{\dsd}$,
where $\dsd$ is the number of
vertices in the graphs $D_1$, $D_2$ which have different cyclic orders. In the
future we will sometimes denote $\tcDbd$ simply as $\tcD$, assuming that the
choice of cyclic order for every 3-valent graph was somehow fixed.

\begin{lemma}
There is a natural isomorphism $\tcDbd \cong \tcDo/\tcDASo$.
\end{lemma}
\proof
Since $\SgroD\subset \SgrD$, then $\Hsx{D}\subset\Hso{D}$. As a result,
$\tcDbd$ may be considered a subspace of $\tcDo$ and thus we have a map
$f\,:\;\tcDbd\longrightarrow \tcDo/\tcDASo$. On the other hand, one can
construct a natural map $g\,:\;\tcDo\longrightarrow \tcDbd$ in the following
way: if a 3-valent graph $D_1$ is isomorphic to a graph $D_2\in\bD$, then $g$
maps $\Hso{D_1}$ to $\Hsx{D_2}\subset\Hso{D_1}$ as
$(-1)^{\dsd}\Pgix{D_1}$. Obviously, $\tcDASo\subset \ker\,g$, so we have a map
$h\,:\;\tcDo/\tcDASo \longrightarrow \tcDbd$. We leave it to the reader to
check that $f$ and $h$ constitute an isomorphism.\qed

After constructing an isomorphism
$h\,:\;\tcDo/\tcDASo\longrightarrow \tcDbd$ we define the space $\tcDIHXx$
simply as the image of
$\tcDIHXo/ (\tcDIHXo\cap\tcDASo)$. Thus we proved the following
\begin{theorem}
\label{t2.1}
There is a natural isomorphism $\cD \cong \tcDbd/\tcDIHXx$.
\end{theorem}
The grading subspaces
$\cD_{m,n}$ turn out to be the quotients $\tcD_{m,n}/(\tcD_{m,n}\cup\tcDIHXx)$.

The advantage of this description of $\cD$ is that it allows us to
work with rather natural spaces $\Invgdx{\SmHor{D}}$ instead of bigger
and less symmetric spaces $\Invgdo{\SmHor{D}}$.

\nsection{Isomorphism between $\cB$ and $\cD$}
\label{s3}

\begin{theorem}
\label{t3.1}
There exists a canonical isomorphism of algebras
\qq
\hxA:\;\cB\longrightarrow \cD,
\label{3.1}
\qqq
which respects the grading
\qq
\hxA:\;\cB_{m,n}\longrightarrow \cD_{m,n-m}.
\label{3.2}
\qqq
\end{theorem}
\begin{corollary}
If $m>n$, then $\cB_{m,n} = \empt$.
\end{corollary}

Before we prove this theorem, we have to establish some facts
concerning the structure of the space $\cB$. We call an edge of a
(1,3)-valent graph \emph{a leg} if this edge is connected to a
1-valent vertex. All other edges are called \emph{internal}.

\begin{lemma}
\label{l3.1}
If two legs of a (1,3)-valent graph $D$ are attached to the same 3-valent
vertex, then $D\in\tcBas$.
\end{lemma}

\proof
Suppose that a (1,3)-valent graph $D$ contains such a 3-valent
vertex. Since the 1-valent vertices of our graphs are not ordered in
any way, then changing the cyclic order at that 3-valent vertex does
not change the graph. Therefore $2D\in\tcBas$ and this proves the
lemma.\qed

Let us call a (1,3)-valent graph \emph{restricted} if each of its
3-valent vertices contains at most one leg. Let $\tcBr$ be a formal
space whose basis vectors are restricted graphs. We introduce
familiar subspaces. The subspaces $\tcBasrii\subset \tcBr$, $i=0,1$
are spanned by the sums of restricted diagrams $D_1,D_2$ which differ
in the ordering at a 3-valent vertex which is attached to $i$ legs.
The subspaces $\tcBihxii\subset \tcBr$, $i=0,1$ are spanned by the
images of the map\rx{1.10} acting on the (3,4)-valent diagrams
whose
single
4-valent vertex contains $i$ legs. Then Lemma\rw{l3.1} has a simple
corollary:
\qq
\cB_{m,n} = \tcBrz_{m,n} / (\tcBasrz + \tcBihxz),\qquad
\mbox{where}\;\; \tcBrz_{m,n} = \tcBr_{m,n}/ (\tcBasro + \tcBihxo) .
\label{3.3}
\qqq
Indeed, this relation follows from the fact that if the 4-valent
vertex of a (3,4)-valent graph $D$ has at least two legs, then the
intersection of the image of the corresponding operator\rx{1.10} with
the space $\tcBr$ is trivial.
Also, it is easy to see that $\tcBasro$ and $\tcBihxo$ are ideals in $\tcBr$,
so the quotient
\qq
\tcBrz = \bopmnzi \tcBrz_{m,n} = \tcBr / (\tcBasro + \tcBihxo)
\qqq
has a graded algebra structure.

Now we begin to construct the isomorphism. Let $D$ be a 3-valent
graph with $N$ edges and cyclic order at vertices.
Thinking of $D$ as a $CW$-complex, let $\Co$ be the space of 1-chains. In other
words, $\Co$ is an $N$-dimensional vector space spanned by the oriented edges
of $D$, if we assume that an edge with the opposite orientation is equal to the
opposite of the edge as an element of $\Co$. Thus, if we pick an orientation on
the edges of $D$, then $\Co$ has a natural basis $e_j$, $1\leq j\leq N$ of the
edges of $D$. We will also need the dual space $\Cos$ with the dual basis
$f_j$, $1\leq j\leq N$. The symmetry group of the graph $G_D$ acts on both
spaces $\Co$ and $\Cos$.


Next, consider a vector space whose basis is formed by
$m$-legged (1,3)-valent restricted graphs such that if we remove
their legs, then we get the 3-valent graph
$D$. We denote the quotient of this space by
its intersection with $\tcBasro$ as $\tcBmd$. We also have to
consider a bigger space. Suppose that we index the edges of $D$ and
then attach $m$ legs to its edges in order to produce restricted graphs. These
(1,3)-valent
graphs still carry the indexing of the edges of $D$. If we factor
this space by its intersection with the obvious analog of $\tcBasro$,
then we get the space $\cchBmdp$. The symmetry group $\SgroD$ of the
graph $D$ acts on $\cchBmdp$ by mapping the edges of $D$ together with their
legs, while preserving the cyclic order at the vertices.
The invariant subspace of this action
is canonically isomorphic to $\tcBmd$:
\qq
\tcBmd = \Invgdo{\cchBmdp}.
\label{3.4*}
\qqq

Let us introduce a
multi-index notation
\qq
\um = (m_1,\ldots,m_N), \qquad |\um| = \sjoN m_j.
\label{3.4}
\qqq
For $N$ non-negative numbers $\um$ and for a choice of orientation of the
edges of $D$
construct a diagram $D_\um$ in
the following way: for every $j$, $1\leq j\leq N$ attach $m_j$ legs
to $D$ on the left side of the edge $e_j$ (the notion of the left
side is well-defined since $e_j$ is oriented). It is easy to see that
all graphs $D_\um$, $|\um|=m$ form a basis of the space $\cchBmdp$,
because after we took the quotient over the analog of the space
$\tcBasro$, we can flip the legs of the graphs of $\cchBmdp$ to a
particular side of each edge of $D$ (at the cost of changing the
signs of the corresponding vectors of $\cchBmdp$).

There is a natural isomorphism
$\xA:\;\cchBmdp \longrightarrow \SmCos$ which acts on the basis
vectors as
\qq
\hxA:\;D_\um\mapsto\pjoN f_j^{m_j}.
\label{3.5}
\qqq

Suppose that the 3-valent graph $D$ has $\Nz$ vertices $v_j$,
$1\leq j\leq \Nz$. Consider the $\Nz$-dimensional space $\Cz$ of 0-chains whose
basis vectors are in a one-to-one correspondence with these vertices.
Then there is a natural boundary map $\xdel:\;\Co\longrightarrow\Cz$.
Let $\chCos$ be the space of 1-cocycles, it is
the subspace of $\Cos$ whose elements annihilate the
kernel of $\xdel$. Apparently,
\qq
\Hor{D} = \Cos/\chCos.
\label{3.6}
\qqq

Let $\cchBihxrodm = \cchBmdp\cap\cchBihxro$, where the space
$\cchBihxro$
is the analog of the space $\tcBihxro$ for the graphs which come from
3-valent graphs with indexed edges.

\begin{lemma}
\label{l3.2}
The map $\hxA$ establishes an isomorphism between the spaces
$\cchBihxrodo$ and $\chCos$.
\end{lemma}

\proof
For $1\leq j\leq \Nz$, denote as $V_j$ the image in $\cchBodp$ of
the operator\rx{1.10} associated with the vertex $v_j$ of $D$ (that
is, one of the two 3-valent vertices in each of the graphs of \fg{f2}
is $v_j$, while the other vertex is attached to a leg). Then the
space $\cchBihxrodo$ is spanned by all the spaces $V_j$.

For $1\leq j\leq \Nz$ and for $x\in\Co$ let $\xdel_j(x)$ be the
coefficient in front of $v_j\in\Cz$ in the expansion of $\xdel(x)$
with respect to the basis $v$.
Then $\ker\xdel = \bigcap_{j=1}^{\Nz} \ker\xdel_j$
and, as a result, the
space $\chCos$ is spanned by the spaces
$V\p_j\subset\Cos$ which annihilate the spaces
$\ker\xdel_j\subset\Co$. It is very easy to see that for every
$j$, $\hxA$ establishes an isomorphism between the corresponding
spaces $V_j$ and $V\p_j$. This proves the lemma.\qed

\begin{lemma}
\label{l3.3}
$\hxA$ establishes the isomorphism between the spaces
$\cchBmdp/\cchBihxrodm$ and
\\
$\SmHord$.
\end{lemma}

To prove this lemma we need a simple fact from linear algebra.
\begin{lemma}
\label{l3.4}
Let $V$ be a finite dimensional vector space and $W$ be its subspace.
Denote by $\Ps$ a symmetrizing projector
$\Ps:\;V^{\otimes m}\longrightarrow S^mV$. Then
\qq
S^m V / \Ps(S^{m-1}V\otimes W) = S^m(V/W).
\label{3.7}
\qqq
\end{lemma}
\proof
We leave the proof to the reader.

\pr{Lemma}{l3.3}
It is easy to see that $\hxA$ maps the space $\cchBihxrodm$ onto
$\Ps(S^{m-1}\Cos\otimes \chCos)$. Then the claim of the lemma follows
from \eex{3.6} and\rx{3.7} if we set $V=\Cos$ and $W=\chCos$ in the
latter equation.\qed

Consider a space $\tcBihxrodm = \tcBmd \cap \tcBihxro$.

\begin{lemma}
\label{l3.5}
There is a natural isomorphism between the quotient spaces
\qq
\tcBmd / \tcBihxrodm = \Invgdo{\cchBmdp /\cchBihxrodm}.
\label{3.8}
\qqq
\end{lemma}

In order to prove this isomorphism we need another linear algebra
lemma.

\begin{lemma}
\label{l3.6}
Let $V$ be a finite-dimensional representation of a finite group $G$.
Let $W\subset V$ be a subspace, which is invariant under the action
of $G$. Then there is a natural isomorphism
\qq
\Invg{V} / \Invg{W}  = \Invg{V/W}.
\label{3.9}
\qqq
\end{lemma}
\proof
For example, one could use the fact that a finite-dimensional
representation of $G$ is a sum of irreducible representations. We
leave the details to the reader.\qed

\pr{lemma}{l3.5}
The cyclic order preserving
symmetry group $\SgroD$ of the 3-valent graph $D$ acts on the space
$\cchBmdp$. Obviously, the symmetrization over this action projects
$\cchBmdp$ onto $\tcBmd$. Thus
\qq
\tcBmd = \Invgdo{\cchBmdp}.
\label{3.10}
\qqq
At the same time, the subspace $\cchBihxrodm$ is invariant under the
action of $\SgroD$ and
\qq
\tcBihxrodm = \Invgdo{\cchBihxrodm}.
\label{3.11}
\qqq
Then \ex{3.8} follows from \ex{3.9} in view of the relations\rx{3.10}
and\rx{3.11}.\qed

Let us introduce a notation $\cBmd = \tcBmd / \tcBihxrodm$.

\begin{corollary}
\label{c3.2}
The map $\hxA$ establishes the isomorphism between the spaces $\cBmd$
and $\Hsio{m}{D}$ (see \ex{1.13} for the definition of the latter
space).
\end{corollary}

\proof This isomorphism follows from the combination of
Lemmas\rw{l3.3} and\rw{l3.5}.\qed

We leave it for the reader to check that the isomorphism $\hxA$ intertwines the
maps
\qq
\cB_{m_1}(D_1)\otimes\cB_{m_2}(D_2) &\longrightarrow &
\cB_{m_1 + m_2} (D_1\cup D_2)
\nonumber\\
\Hsio{m_1}{D_1}\otimes\Hsio{m_2}{D_2} &\longrightarrow &
\Hsio{m_1 + m_2}{D_1\cup D_2}
\qqq
which come from the multiplications in the algebras $\tcB$ and $\tcD$ as
defined in Section\rw{gs}.


\pr{Theorem}{t3.1}
According the definition\rx{3.3} of the space $\tcBrz_{m,n}$,
\qq
\tcBrz_{m,n+m} = \bigoplus_{D:\;\chi(D)=n} \cBmd,
\label{3.12}
\qqq
while by its definition
\qq
\tcDo_{m,n} = \bigoplus_{D:\;\chi(D)=n} \Hsio{m}{D}.
\label{3.13}
\qqq
It is easy to see that $\hxA$ establishes the isomorphisms
\qq
\hxA:\; \tcBasrz\cap\tcBrz_{m,n+m}\longrightarrow
\tcDASo\cap\tcDo_{m,n},\qquad
\tcBihxz\cap\tcBrz_{m,n+m}\longrightarrow\tcDIHXo\cap\tcDo_{m,n}.
\label{3.14}
\qqq
Then \ex{3.2} follows from \eex{3.3} and\rx{1.24*} together with the
isomorphism of Corollary\rw{c3.2}.\qed

\nsection{Rationality Conjecture}
\label{s4}

Recall that Kontsevich integral of a knot $\cK\in S^3$ is a sequence
of vectors $\IB_{m,n}(\cK)\in \cB_{m,n}$, $m\geq 0$, $n\geq m$
depending on the topological class of $\cK$. The space $\cB_{0,0}$ is
1-dimensional, its basis vector is the empty graph, so it can be
naturally identified with $\IC$. It is known that
$\IB_{0,0}(\cK) = 1$.

We combine the vectors $\IB_n(\cK)$ into a formal
power series of a formal variable $\hb$
\qq
\IB(\cK;\hb) = 1 + \smnoge \IB_{m,n}(\cK)\,\hb^n\in\cB.
\label{4.1}
\qqq
Prior to formulating a conjecture about the structure of
$\IB(\cK;\hb)$
we have to apply to it some transformations. First, we apply
the wheeling map $\hOm:\;\cB\longrightarrow\cB$, described in\cx{Wh},
in order to produce
\qq
\IO(\cK;\hb) = \hOm(\, \IB(\cK;\hb)) =
1 + \smnog \IO_{m,n}(\cK)\,\hb^{m+n} \in \cD,\qquad
\IO_{m,n}\in\cB_{m,n}.
\label{4.1*}
\qqq
Then we apply the isomorphism $\hxA$, which maps Kontsevich integral
from $\cB$ to $\cD$. More precisely, we choose a set
$\bD$ of 3-valent graphs $D$ such
that each type of a graph (without distinguishing them by cyclic
order at vertices) is represented there exactly once, and then we map $\cB$ to
$\cD_{\bD}$ as described at the end of Section\rw{gs}.
Thus we get
\qq
\ID(\cK;\hb) = \hxA(\, \IB(\cK;\hb)) =
1 + \smnog \ID_{m,n}(\cK)\,\hb^{m+n} \in \cD,
\qquad \ID_{m,n}(\cK)\in \cD_{m,n}.
\label{4.2}
\qqq
%


By using the algebra structure of $\cD$ and manipulating the formal power
series in $\hb$ we can define the logarithm of Kontsevich integral
\qq
&
\Il(\cK;\hb) = \log \ID(\cK;\hb) =
\smnog \Il_{m,n}(\cK)\,\hb^{m+n} \in \cD,
\qquad \Il_{m,n}(\cK)\in \cD_{m,n}.
\label{4.3*}
\qqq
through the formula
\qq
\log(1+x) = \sum_{n=1}^{\infty} (-1)^{n-1}{x^n\over n}
\qqq
The advantage of the logarithm $\Il(\cK;\hb)$ is that it can be expressed
exclusively in terms of \emph{connected} 3-valent graphs.

Kontsevich integral $\Il(\cK;\hb)$ belongs to the quotient space\rx{1.23}. Let
$\tIl(\cK;\hb)$ be a representative of $\Il(\cK;\hb)$ in the space
$\tcD$ (Of course,
it is defined only up to an element
of $\tcDIHX$).
We present
$\tIl(\cK;\hb)$ as
\qq
\tIl(\cK;\hb) = \sDc \smzi x_m(\cK,D)\, \hb^{\chi(D)+m},
\label{4.4}
\qqq
where $\bDc\subset \bD$ is a subset of connected 3-valent graphs and
$x_m(\cK,D)\in \Hsm{D}$.

Now we are almost ready to formulate our conjecture. Let $V$ be a vector space.
For $x\in V$ we define $e^x\in \Ss V$ by the power series $e^x =\snzi x^n/n!$.
If $\Lambda$ is a lattice in $V$, then we extend this exponential map to
an injection of a group algebra $\Exp:\;\IQ[\Lambda]\rightarrow \Ss V$. For a
graph $D$, $\Hoz{D}$ forms a lattice in $\Hor{D}$.
We denote
\qq
\HsQ{D} = \Exp
\Qalghz
\subset\Hs{D}.
\qqq
In other words, $\HsQ{D}$ is $\SgrD$-invariant part of the rational span of the
exponents of the elements of $\Hoz{D}$ and $\Exp$ establishes its isomorphism
with $\Qalghz$.

Now recall that if $D$ has $N$ edges, then $e_j$ ($1\leq j\leq N$)
denote the oriented edges forming a basis in the space of 1-chains $\Co$, while
$f_j$, $1\leq j\leq N$ form the dual basis in the dual space $\Cos$.
In view of \ex{3.6} we can think of $f_j$ as elements of $\Hor{D}$.

\begin{lemma}
The product of the Alexander-Conway polynomial of $e^{f_j}$ is an element of
the algebra
$\HsQ{D}$:
\qq
\pjoN \APbas{\cK}{\exp(f_j)}\in\HsQ{D}.
\label{4.4*}
\qqq
and its inverse is a well-defined element of $\Hs{D}$
\end{lemma}
\proof
To prove relation\rx{4.4*}, we have to show that its \lhs is $G_D$-invariant.
The elements of the group $G_D$ not only permute $f_j$, $1\leq j\leq N$, but
they may also reverse the orientation of some edges of $D$ and thus change the
signs of corresponding $f_j$. However, the relation
\qq
\APbas{\cK}{1/t} = \APbas{\cK}{t},
\label{4.5*}
\qqq
guarantees that this change of sign does not affect the expression\rx{4.4*}
and hence it is $G_D$-invariant.
At the same time, the Alexander-Conway polynomial satisfies the
property $\APbas{\cK}{1}=1$ which guarantees that the inverse of\rx{4.4*} can
be inverted within $\Hs{D}$.\qed


Let us introduce a notation
\qq
\Il(\cK,D) = \smzi x_m(\cK,D) \in \Invgd{\SHor{D}}.
\label{4.5*1}
\qqq
The only 3-valent graph $D$ with $\chi(D)=0$ is a
 circle.
The value
of $\Il(\cK,\crcl)$ has been established by D.~Bar-Natan and
S.~Garoufalidis in\cx{BG}
\qq
\Il(\cK,\crcl)
= \hlfv
\lrbs{
\log \lrbc{\sinh(f/2) \over (f/2) } -
\log \APbas{\cK}{\exp(f)}
},
\label{4.6}
\qqq
where $f$ represents the integral generator of $\Hor{\crcl}$. Our conjecture
deals with the value of $\Il(\cK,D)$ for graphs with $\chi(D)\geq 1$. Recall
that such graphs have exactly $N=3\chi(D)$ edges.

\begin{conjecture}
\label{c4.1}
The representative $\tIl(\cK;\hb)\in\tcD$ of Kontsevich integral
$\Il(\cK;\hb)\in\cD$ can be chosen in such a way that for any
$D\in\bD$, $\chi(D)\geq 1$ there exists an element
$y(\cK,D)\in\HsQ{D}$
such that
\qq
\Il(\cK,D)
= {y(\cK,D) \over \pjochD \APbas{\cK}{\exp(f_j)} }.
\label{4.5}
\qqq
\end{conjecture}
\begin{remark}
\rm
Andrew Kricker has proved this conjecture in his paper\cx{Kr}.
\end{remark}

\begin{remark}
\rm
D.~Thurston presented arguments which show that if Conjecture\rw{c4.1} is true
as it is formulated, then it should also be true if one defines $\ID(\cK,\hb)$ directly
as an image of $\IB(\cK,\hb)$ under the isomorphism $\hxA$ without appying the wheeling
map $\hOm$ of \ex{4.1*}.
\end{remark}

\begin{remark}
\rm
It is convenient to introduce some other notations in relation to \ex{4.5}. Let
$p(\cK,D)\in \Qalghz$ be such that $\Expalg{p(\cK,D)} = y(\cK,D)$. Also, if
we index the edges of $D$ in such a way that $\fbs$ form a
basis of
$\Hoz{D}$ and $\Hor{D}$, then we can write $\Il(\cK,D)$ and
$y(\cK,D)$ more explicitly as
\qq
\Il(\cK,D) & = & \Il(\cK,D;\fbs),
\label{4.5*2}\\
p(\cK,D) & = & p(\cK,D;\fbs),
\label{4.5*2*1}\\
y(\cK,D) & = & p(\cK,D;\efbs),
\label{4.5*2*2}
\qqq
where
\qq
&\Il(\cK,D;\xbs)\in\IQ[[\xbs]],
\nonumber\\
& p(\cK,D;\tbs)\in
\IQtbs.
\label{4.5*3}
\qqq

\end{remark}

\nsection{Rational structure of the Jones polynomial}
\label{s5}

There is a well-known relation between the Kontsevich integral and the colored
Jones polynomial of a knot, so the rationality conjecture\rw{c4.1} should
manifests itself in the structure of the latter object. In fact, this
manifestation observed in\cx{Ro}, served for us as evidence which led to the
rationality conjecture. Another advantage in establishing a relation between
\ex{4.5} and the rational expansion of the Jones polynomial\cx{Ro} is that
at present it is much easier to calculate the colored Jones polynomial than
Kontsevich integral. Therefore, working out the rational expansion of\cx{Ro} is
a practical way of finding the polynomials $y(\cK,D)$ of \ex{4.5}.

Let us recall the exact relation between the Kontsevich integral and a colored
Jones (or, more generally, HOMFLY) polynomial based on a simple Lie algebra
$\gg$. We equip $\gg$ with the ad-invariant
%
scalar product
normalized in such a way that long roots have length $\sqrt{2}$ (this
scalar product allows us to identify the dual space $\gg^*$ with
$\gg$ itself). Let
$\val\in\gh$ be the hightest weight of a representation of $\gg$,
shifted by $\vr$ (which is half the sum of positive roots of $\gg$).
Reshetikhin and Turaev associate to this data a polynomial
$\JvaK\in\ZZ[q^{\pm 1/2}]$. If we substitute
\qq
q = e^{\hb},
\label{5.1*}
\qqq
then we
can expand $\JvaK$ in power series of $\hb$
\qq
\JvaK = \snzi \pnvaK\,\hb^n,
\label{5.1}
\qqq
whose coefficients $\pnvaK$ are polynomials of $\val$.
The same series\rx{5.1} can be deduced from the value of Kontsevich
integral.

The data $\gg,\val$ defines an element in the dual space
$\cBs$, which is called \emph{the weight system}. We will define it
in such a way that it will be
suitable for application to
$\IO(\cK;\hb)$.
The first steps in the definition of the weight systems are fairly
standard. Let $\vx_a$, $1\leq a\leq \dim \gg$ be a basis of $\gg$.
Define the structure constants $f_{abc}$ by the relation
\qq
[\vx_a,\vx_b] = \sum_{c=1}^{\dim \gg} \fabco\, \vx_c.
\label{5.2}
\qqq
We can raise and lower the indices of $\fabco$ with the help of the
metric tensor
\qq
h_{ab} = \vx_a\cdot \vx_b
\label{5.3}
\qqq
and its inverse $h^{ab}$.

Let $D$ be a (1,3)-valent graph, $\dego(D)=m$, $\degt(D)=n+m$.
Suppose that if we strip off its legs, then we get a 3-valent graph
$\Dz$. Let us orient the edges of $\Dz$ and assign
orientation to the edges of $D$ in such a way that it is compatible
with the orientation of $\Dz$ and legs are oriented in the direction
from 1-valent vertex to 3-valent vertex. Next, we assign the tensors
$f$ to 3-valent vertices, assigning their indices to attached edges
according to the cyclic ordering. We use the upper indices for the
incoming edges and lower indices for the outgoing edges. Finally, we
take the product of all tensors $f$ assigned to 3-valent vertices,
contract each pair of indices of $f$'s along each internal edge,
while contracting each index assigned to a leg with $\a_a$
($\val = \sum_{a=1}^{\dim\gg} \a_a\,\vx_a$).
Thus we get a Weyl group invariant
homogeneous polynomial $\pDval$ of $\val$ of degree $m$. It is easy
to see that it does not depend on the choice of orientation of the
edges of $\Dz$. For a fixed weight $\val$, $\yp$ assignes a number to
each (1,3)-valent graph, so $\yp\in\tcB^*$. In fact, due to the
anti-symmetry of $f$ and to the Jacobi identity, satisfied by
the commutator\rx{5.2}, $\yp$ annihilates the subspaces $\tcBas$ and
$\tcBihx$ and therefore it can be projected to $\cBs$.

The usual way to proceed further is to convert $\pD$ as a Weyl group invariant
polynomial on $\gh$ into an element of $\Inv{\gg}{S^m\gg}$, then use a PBW map
to convert it into an element of $\Inv{\gg}{U\gg}$ and calculate the trace of
that element in a $\gg$ module with the highest weight $\val - \vr$, thus
obtaining another polynomial $\pval$ of $\val\in\gh$ which is the standard
weight of the graph $D$ coming from $\gh$, or thinking of it as a function on
all graphs $D$, $\yps$ is a weight system on $\cB$. Then the relation between
the expansion\rx{5.1} and Kontsevich integral is
\qq
\JvaK = \dval \lrbc{
1 + \smnog \yps(\,\IB_{m,n}(\cK),\val)\,\hb^{m+n} },
\label{5.4*}
\qqq
where $\dval$ is the dimension of the representation of $\gg$ with
the shifted highest weight $\val$. However, as explained in\cx{Wh}, the
wheeling map allows one to get the expansion\rx{5.1} straight from the weight
$\pDval$ without going through PBW map and calculating the trace:
\qq
\JvaK = \dval \lrbc{
1 + \smnog \yp(\,\IO_{m,n}(\cK),\val)\,\hb^{m+n} },
\label{5.4}
\qqq
This is the formula that we will work with, because the weight function
$\pDval$ is easy to transfer from $\cB$ to $\cD$.
The inverse of the dual isomorphism map $\hxA^*$ maps the weight
system $\yp\in \cBs$ into an element of $\cDs$, which we will call
$\ypD$. In order to see how $\ypD$ acts on $\cD$ we come back to the
calculation of $\pDval$ and modify it.

Suppose that $\gg$ has $2k$ roots $\l_1,\ldots,\l_{2k}$. Let us index them in
such a way that $\l_1,\ldots,\l_{\kg}$ are positive roots and
$\l_1,\ldots\l_{\rg}$ are simple roots, $\rg$ being the rank of $\gg$.

For a root $\l$ of $\gg$ let $\Prl$
denote the operator projecting $\gg$ onto the root space
$V_\l\subset \gg$. We also introduce an operator $\Ph$, projecting
$\gg$ onto $\gh$. Let us assign a root of $\gg$
or the Cartan subalgebra to each internal edge of $D$. Let $\tbfS$ be
a set of all such assignments. For an assignment $c\in\tbfS$ we
modify the contraction of indices of tensors $f$ in the following
way: if an internal edge carries an index $a$ at the beginning and
index $b$ at the end, then instead of contracting them (that is,
instead of setting $a=b$ and taking a sum over their values) we bring
in an extra factor $P_b^a$, where $P$ is the projector corresponding
to the subspace assigned to that edge by $c$, and then contract the
pairs of indices $a$ and $b$ independently. In other words, we
project Lie algebras $\gg$ flowing along the internal edges of $D$
onto root spaces and Cartan subalgebras. Let us denote the resulting
number as $\pcDval$. Since the sum of projectors $\Ph$ and $\Prl$ for
all roots $\l$ of $\gg$ is equal to the identity operator, then
\qq
\pDval = \sum_{c\in\tbfS} \pcDval.
\label{5.5}
\qqq

The sum in the \rhs of this equation can be simplified.
Since $\val\in \gh$, then
\qq
[\val,\vy] = (\val\cdot\l)\,\vy\qquad\mbox{if $\vy\in V_\l$},
\qquad [\val,\vy] = 0 \qquad\mbox{if $\vy\in\gh$}.
\label{5.6}
\qqq
Therefore, $\pcDval=0$ unless the following two conditions are
met. First, $c$ must assign the same projector to internal edges of
$D$ which correspond to the same edge of $\Dz$. Second, there is a
\emph{compatibility requirement} at every 3-valent vertex: Cartan
subalgebra can be assigned to at most one of its edges and the sum of
the roots on incoming edges is equal to the sum of the roots on
outgoing edges.
Thus we can replace
the set $\tbfS$ in \ex{5.5} with the set $\bfS$ of `compatible' assignments
whose elements assign
subspaces to the edges of $\Dz$ in such a way that the compatibility
condition is satisfied at all of its vertices.

Equations\rx{5.6} also indicate that the effect of leg contractions
is easy to take into account in the calculation of $\pcDval$. If a leg is
attached to at least one edge, to which a Cartan subalgebra is
assigned, then $\pcDval=0$. Otherwise, if $m_j$ legs are attached on
the left side of an oriented edge $e_j$ of $\Dz$ to which a root
$\l$ is assigned, then they contribute a factor of
$(\val\cdot\l)^{m_j}$. Let $\lcj$ denote the root of $\gg$ assigned
by $c\in\bfS$ to the edge $e_j$ of $\Dz$. If $c$ assigns $\gh$ to
$e_j$, then we set $\lcj=0$. With these notations we see that
\qq
\pcDval = \pcDz \pjoN (\val\cdot\lcj)^{m_j},
\label{5.7}
\qqq
where $\pcDz = \pcDzval$ (we had to introduce this new notation
because the graph $\Dz$ has no legs and as a result $\pcDzval$ does
not depend on $\val$). Note that in \ex{5.7} we adopted a convention
that $0^0=1$.

The isomorphism\rx{3.5} completes the translation of $\pcDval$ into
the language of 3-valent graphs. For an assignment $c\in\bfS$
consider a linear combination of edges
\qq
\ecval = \sjoN (\val\cdot\lcj)\,e_j\in \Co.
\label{5.8}
\qqq
According to the compatibility condition satisfied by $c$,
$\ecval\in\ker\xdel=\oHor{D}$. Therefore, we can evaluate an element
$x\in \SsHor{D}$ on $\ecval$ and get a number (or a formal series)
$x(\ecval)$. Equations\rx{3.5},\rx{5.7} and\rx{5.8} indicate that for
an
element $x\in\cchBmdpz/\tcBihxri{1}(D_0,m)$,
\qq
\pcxval = \pcDz\,(\hxA\,x)(\ecval).
\label{5.9}
\qqq
Then, according to \ex{5.5}, after taking a sum over the assignments
of $\bfS$, we come to the following relation: for any $x\in\tcBmdz$,
\qq
\pxval = \pDxval,
\label{5.10}
\qqq
where
\qq
\pDyval=\sum_{c\in\bfS}\pcDz\, y(\ecval),
\qquad y \in \Hsio{m}{D_0}.
\label{5.11}
\qqq
Thus \ex{5.11} defines the element $\ypD\in\cDs$ corresponding to
$\yp\in\cBs$.

Applying \ex{5.10} to \ex{5.4}, we find that
\qq
\JvaK = \dval \lrbc{
1 + \smnog \ypD(\,\ID_{m,n}(\cK),\val)\,\hb^{m+n} }.
\label{5.12}
\qqq
It is easy to see that the weight system $\ypD$ behaves nicely under
the multiplication of elements of $\cD$:
$\ypD(xy,\val) = \ypD(x,\val)\,\ypD(y,\val)$ for any $x,y\in\cD$.
Therefore the analog of \ex{5.12} holds for the modified
integral\rx{4.3*}
\qq
\lgdval
= \smnog \ypD(\,\Il_{m,n}(\cK),\val)\,\hb^{m+n},
\label{5.13}
\qqq
and for its representative\rx{4.4} in the space $\tcD$
\qq
\lgdval = \sD \smzi \ypD(\,x_m(\cK,D),\val)\, \hb^{\chi(D)+m}.
\label{5.14}
\qqq
%
By using the formula\rx{5.11} for the
weight system, we can rewrite \ex{5.14} as
\qq
\lgdval = \sD \sum_{c\in\bfS}
\smzi \pcD\,x_m(\cK,D)(e_{c,\val})\, \hb^{\chi(D)+m},
\label{5.15}
\qqq
where $x_m(\cK,D)(e_{c,\val})$ denotes the evaluation of the element
$x_m(\cK,D)\in \Hsm{D}$ on $e_{c,\val}\in \oHor{D}$. According to
\ex{5.8}, $e_{c,\val}$ is a linear function of $\val$, while
$x_m(\cK,D)(e_{c,\val})$ is the homogeneous polynomial of
$e_{c,\val}$
of degree $m$. Therefore, \ex{5.15} can be further modified as
\qq
\lgdval & = & \sD \hb^{\chi(D)}\sass
\pcD\smzi x_m(\cK,D)(e_{c,\hb\val})
\nonumber\\
& = &
\sD \hb^{\chi(D)}\sum_{c\in\bfS}
\pcD\Il(\cK,D) (e_{c,\hb\val})
\label{5.16}
\qqq
the last line coming from \ex{4.5*1}.
Since by the definition of the
dual basis $f_j(e_i) = \delta_{ij}$, then according to \ex{5.8},
$f_j(e_{c,\hb\val}) = \hb\,(\val\cdot\lcj)$ and as a result, in view
of\rx{5.1*},
\qq
\Il(\cK,D) (e_{c,\hb\val}) =
\IlKi{D}{\aockh}
\label{5*.1}
\qqq
(see \ex{4.5*2}). Equation\rx{4.6} allows us to write the contribution of the
`1-loop' graph ($\chi(D)=0$) explicitly. Assignments
$c$ simply put different roots
on the circle,
$\pccir=1$ and
\qq
\lefteqn{
\sass \pccir \IlKi{\crcl}{\hb(\val\cdot\lci{1})}
}
\label{5*.2}
\\
&&
\hspace{1in}
=
\sum_{j=1}^k \log \lrbc{ q^{(\val\cdot\l_j)/2} -
q^{-(\val\cdot\l_j)/2}
\over \hb \,(\val\cdot\l_j)}
- \sum_{j=1}^k \log \APbas{\cK}{q^{\val\cdot\l_j}}.
\nonumber
\qqq
Thus if we exponentiate both sides of \ex{5.16} and use the
formulas\rx{5*.1},\rx{5*.2} and the dimension formula
\qq
d_{\val} = \pjok { \val\cdot \l_j \over \vr \cdot \l_j},
\qqq
then we find that
\qq
\JvaK  =
{\dqval \over \Dg(\cK;\qaok)}\;\;\Cqg \;
\exp\lrbc{\snoi \JlKn\,\hb^n }
\label{5*.3}
\qqq
where
\qq
&\JlKn = \sDln \;\;\sass\pcD\,\IlKi{D}{\aockh},
\label{5*.3*1}
\qqq
while
\qq
&\dqval = \pjok {q^{(\val\cdot\l_j/2)} - q^{-(\val\cdot\l_j)/2}
\over q^{(\vr\cdot\l_j/2)} - q^{-(\vr\cdot\l_j)/2} }
\label{5.20*2}
\qqq
is called the quantum dimension of the $\gg$-module with highest weight
$\val-\vr$ and
\qq
&\Cqg = \sum_{j=1}^k \log \lrbc{ q^{(\vr\cdot\l_j)/2} -
q^{-(\vr\cdot\l_j)/2}
\over \hb \,(\val\cdot\l_j)} = 1 + \cO(\hb^2),
\label{5.20*4}\\
&\Dg(\cK;\qaok) = \pjok \APbas{\cK}{q^{\val\cdot\l_k}}.
\label{5.20*5}
\qqq

Now let us apply Conjecture\rw{c4.1} to the \rhs of \ex{5*.3*1}. According to
\eex{4.5} and\rx{4.5*2},
\qq
\IlKi{D}{\aockh}
=
{
p(\cK,D;q^{\val\cdot\lci{1}},\ldots,q^{\val\cdot\lci{\chi(D)+1}})
\over \pjochD \APbas{\cK}{q^{\val\cdot\lcj}} }
\qqq
Therefore if we bring all terms in the sums of \ex{5*.3*1} to the common
denominator $$\Dg^{3n}(\cK;\qaok),$$ then we find that $\JlKn$ has a rational
form
\qq
\JlKn = {\plgn(\cK;\qaok)\over \Dg^{3n}(\cK;\qaok)},\qquad
\plgn(\cK;\tok)\in\IQtorg.
\label{5*.4}
\qqq
Then substituting this formula to \ex{5*.3}, exponentiating the formal power
series and expanding $\Cqg$ in powers of $\hb$ we come to the following
\begin{crl}
For a knot $\cK$ and a simple algebra $\gg$ there exist the
polynomials
\qq
&p_n(\cK;t_1,\ldots,t_{\rg}) \in \IQtorg,
&\qquad n\geq 0,
\label{5.19*}
\qqq
such that
\qq
\JvaK & = &
{\dqval \over \Dg(\cK;\qaok)}
\lrbc{1+
\snoi   {p_n(\cK;\qaok) \over \Dg^{3n}(\cK;\qaok) }\,\hb^n
}.
\label{5.20*}
\qqq
\end{crl}

We can check this prediction for the case of $\gg=su(2)$. In fact, in this case
the power of $\Dg$ in denominators\rx{5.20*} can be reduced. Indeed, the algebra
$su(2)$ has only one positive root. As a result, the elements of
$\bfS$ assign the subspaces of $su(2)$ to the edges of a graph $D$
in such a way that for any three edges attached to the same vertex,
two are assigned a root space and the third is assigned the Cartan
subalgebra.
Therefore, of $3\chi(D)$ edges that a graph $D$ has, $\chi(D)$ edges always
carry a Cartan subalgebra and only $2\chi(D)$ edges carry the root spaces.
%
Therefore, in case of $su(2)$ \ex{5.20*} is reduced to
\qq
\JaK = {[\a]\over \APKqa}
\lrbc{ 1 +
\snoi
{p_n(\cK;\qa) \over \APKqan }\,\hb^n
},
\label{5.22}
\qqq
where $\a$ is the dimension of the $su(2)$ module attached to the knot $\cK$
and
\qq
[\a] = {q^{\a/2} - q^{-\a/2}\over q^{1/2} - q^{-1/2} }
\label{5.25}
\qqq
is its quantum dimension.


Equation\rx{5.22} can be verified directly. We proved in\cx{Ro} that
for a knot $\cK$ in $S^3$ there exist the polynomials
\qq
\PnKt \in \ZZ[t^{\pm 1}],\;\;n\geq 1,
\label{5.23}
\qqq
such that the expansion\rx{5.1} can be rewritten as
\qq
\JaK = {[\a]\over \APKqa}
\lrbc{1+
\snoi {\PnKqa \over \APKqan}\,h^n,
},
\label{5.24}
\qqq
where
\qq
h = q-1 = e^\hb - 1.
\label{5.25*}
\qqq
It is easy to see that \ex{5.22} follows easily from \ex{5.24}.

\nsection{2-loop invariant and the $SU(3)$ colored Jones polynomial}
\label{s6}

Let us describe more precisely the implications of Conjecture\rw{c4.1} for the
value of Kontsevich integral at the level of `2-loop' graphs, \ie the graphs
with $\chi(D)=1$. There are only 2 such connected graphs in $\bDc$:
the theta-graph $D_1$
and the dumbbell $D_2$ of \fg{f3}.
\begin{figure}[hbt]
\leavevmode \centerline{
\hspace*{1.1in}
\epsfbox{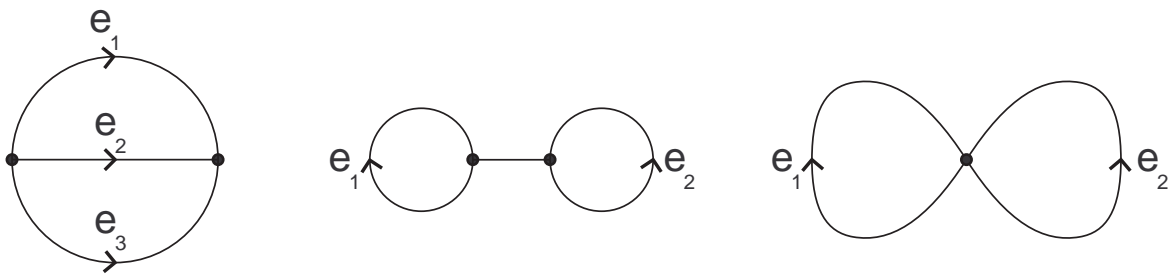}} \caption{The 2-loop
graphs $D_1$, $D_2$ and $D_3$}
\label{f3}
\end{figure}
Therefore, we can present the 2-loop part of the Kontsevich integral\rx{4.3*}
as
\qq
\smzi\Il_{m,1}(\cK) = \Il(\cK,D_1;\fDot{1}) +
\Il(\cK,D_2;\fDot{2})
\label{6.1}
\qqq
(\cf \eex{4.4},\rx{4.5*1} and\rx{4.5*2}), where we used a notation $\fD{i}{j}$
instead of simply $f_i$ in order to distinguish the dual edges coming from
different graphs. The formal power series in the \rhs of \ex{6.1} are not
themselves the invariants of $\cK$. They become the invariants only after the
factorization over the subspace $\tcDIHXx$ (see Theorem\rw{t2.1} and preceding
discussion). Let us describe the IHX indeterminacy in these power series more
precisely. The graph $D_3$ of Fig.\rw{f3} is the only connected 2-loop
graph with a 4-valent vertex. Applying the operator $\hdihx$ of\rx{1.22} to an
element $z(\fDot{3})\in\SHor{D_3}$ we get
\qq
\lefteqn{
{2\over 3}[z(\fDot{1}) + z(\fD{2}{1},-\fD{1}{1}-\fD{2}{1})
+ z(-\fD{1}{1}-\fD{2}{1},\fD{1}{1}) ]
-z(\fDot{2})
}&&
\nonumber\\
&&\hspace{4in}\in\bigoplus_{i=1}^2 \,\Invgdix{i}{\SHor{D_i}}.
\label{6.2}
\qqq
%
In this formula we assumed for simplicity of notation that
$z(x_1,x_2)\in\IQ[[x_1,x_2]]$ already has the symmetries
\qq
z(x_1,x_2) = z(x_2,x_1) = -z(-x_1,x_2),
\label{6.3}
\qqq
which makes the additional symmetrization of the expression\rx{6.2}
unnecessary. Expression\rx{6.2} indicates that by using the IHX freedom we can
bring the expression\rx{6.1} to the form
\qq
&\smzi\Il_{m,1}(\cK) = \zthK{\fDot{1}}\in\Invgdix{D_1}{\SHor{D_1}},
\label{6.4}
\qqq
where
\qq
\lefteqn{
\zthK{x_1,x_2} = \Il(\cK,D_1;x_1,x_2)
}
\label{6.5}\\
&& +
{2\over 3}[
\Il(\cK,D_2;x_1,x_2) + \Il(\cK,D_2;x_2,-x_1-x_2) +
\Il(\cK,D_2;-x_1-x_2,x_1)],
\nonumber
\qqq
thus eliminating the graph $D_2$ from Kontsevich integral. At the same time,
expression\rx{6.2} shows that $\zthK{x_1,x_2}$ of \ex{6.4} is the
$IHX$-invariant combination and therefore it is the only 2-loop invariant of
$\cK$.

The rationality conjecture implies that $\zthK{x_1,x_2}$ also has a
rational structure. Indeed, according to the conjecture, one can use the IHX
freedom in order to bring the terms in the \rhs of \ex{6.1} to the following
form:
\qq
\Il(\cK,D_1;\xot) & = &
{p(\cK,D_1;\exot) \over \APKex{1}\APKex{2}\APKe{x_1+x_2}}
\nonumber\\
\Il(\cK,D_2;\xot) & = & {p(\cK,D_2;\exot) \over \APKex{1}\APKex{2}}.
\label{6.6}
\qqq
Then according to \ex{6.5}, $\zthKx$ has a form
\qq
&\zthKx = {\pthKex\over \APKex{1}\APKex{2}\APKe{x_1+x_2}},
\label{6.7}
\qqq
where
the polynomial $\pthKt\in \IQtot$
is an invariant of $\cK$. Both this polynomial and a rational function
\qq
\zthsKt = {\pthKt\over \APKti{1}\APKti{2}\APK{t_1 t_2} }
\label{6.7*}
\qqq
have the symmetries
%
\qq
f(\tot) = f(t_2,t_1) = f( (t_1 t_2)^{-1},t_2) = f(t_1^{-1},t_2^{-1})
\label{6.8}
\qqq
implied by
the symmetry group $\Sgri{D_1}$. Finally, we rewrite \ex{6.4} with
the help of \ex{6.7}
\qq
\smzi\Il_{m,1}(\cK) =
\zthK{\fDot{1}} =
{\pthKef \over \APK{\efD{1}}\APK{\efD{2}}\APK{e^{-f_{1,D_1}-f_{2,D_1}}}  }.
\label{6.8*}
\qqq

It is easy to see from its definition that Kontsevich integral\rx{4.1} does not
contain (1,3)-valent graphs without legs. The wheeling map $\hOm$ produces such
graphs, however their Euler characteristic is at least 2. Therefore,
$\Il_{0,1}(\cK)=0$ in\rx{4.3*} and this means that
\qq
\zthK{0,0} = \zthsK{1,1} = 0.
\label{6.8*1}
\qqq

The polynomial $\pthKt$ can be extracted from the colored $SU(3)$ Jones
polynomial as described in Section\rw{s5}. In\cx{Rop} we will prove
a slightly strengthened version of \ex{5.20*}
for the groups $SU(n)$:
\qq
&
\JvaK  =
{\dqval\over \Dg(\cK;\qaok)}
\lrbc{1 +
\snoi   {P_n(\cK;\qaok) \over \Dg^{3n}(\cK;\qaok) }\,h^n
},
\label{6.9*}
\qqq
\qq
&
P_n(\cK;\tok) \in \ZZtok
\nonumber
\qqq
(note that here we used an expansion parameter $h=e^{\hb}-1$ instead of $\hb$
and as a result obtained the polynomials with \emph{integer} coefficients). For
the case of $SU(3)$ this formula implies that
\qq
&
\JvaK = {\dqval \over \Dg(\cK;\qaot)}\lrbc{1 + \hb\,
\FoK{\qaot} + \cO(\hb^2)},
\label{6.10*}
\qqq
where
\qq
& \FoK{\tot} =
{P_1(\cK;\tot)\over
\lrbs{
\APKti{1}\APKti{2}\APK{t_1 t_2}
}^3
}.
\label{6.11*}
\qqq
As we explained in Section\rw{s5}, a similar formula\rx{5*.3} can be obtained
by applying the $su(3)$ weight system to the logarithm of the Kontsevich
integral of $\cK$. Comparing \eex{6.10*} and\rx{5*.3}
and taking into account that $\Cqg = 1 + \cO(\hb^2)$,
we see that
\qq
\FoK{\qaot} = \JlKo.
\label{6.12*}
\qqq
Equations\rx{5*.3*1},\rx{6.7*} and\rx{6.8*} show that
\qq
\JlKo = \sass\pcDo \zthsK{\qaotc}
\label{6.13*}
\qqq
and the sum in this formula goes over the compatible assignments of root spaces
and Cartan subalgebra to the egdes of the
$\theta$-shaped graph $D_1$. There are two types of
such assignments. The first
one assigns two opposite roots to two edges and Cartan subalgebra to the third
edge, so
$\pcDo=2$. There are 3 choices of pairs of roots, and within each choice
there are 6 distinct assignments which give the same contributions due to the
symmetries\rx{6.8}. Therefore, the total contribution of the first assignment
to the \rhs of \ex{6.13*}
is
\qq
12
\Big(
\zthsK{\qaic{1},1} + \zthsK{\qaic{2},1} +
\zthK{\qalot,1 }
\Big).
\label{6.14*}
\qqq
Assignments of the second type put 3 different roots on the edges of $D_1$, so
$\pcDo=1$ There are 2 choices of compatible triplets of roots, and there are 6
ways to assign each triplet to the egdes of $D_1$. Thus we have 12 assignments
of the second type, and each gives the same contribution due to the
symmetries\rx{6.8}. Therefore, the total contribution of the second assignment
to the \rhs of \ex{6.13*} is
\qq
12\,\zthsK{\qaot}.
\label{6.15*}
\qqq
Thus putting the sum of\rx{6.14*} and\rx{6.15*} in the \rhs of \ex{6.13*} we
find from \ex{6.12*} that
\qq
\FoK{\tot} =
12\,\Big(
\zthsK{t_1,1} + \zthsK{t_2,1} + \zthsK{(t_1 t_2)^{-1},1} + \zthsKt
\Big).
\label{6.16*}
\qqq
It is easy to solve this equation for $\zthsKt$. By setting $t_2=1$ and using
\ex{6.8*1} and the symmetries\rx{6.8} we get
\qq
\FoK{t_1,1} = 36 \zthsK{t_1,1},
\label{6.17*}
\qqq
hence
\qq
\zthsKt = {1\over 36}\;
\Big(
3\FoK{\tot} - \FoK{t_1,1} - \FoK{t_2,1} - \FoK{(t_1 t_2)^{-1},1}
\Big).
\label{6.18*}
\qqq
In\cx{Rop} we will present a relatively efficient way of calculating
$\FoK{\tot}$. We have already written a Maple V program\cx{RoM}
which implements this
algorithm. For a knot presented as a cyclic closure of a braid, this program
calculates $\APKt$, $P_1(\cK;\tot)$ of \ex{6.11*} and then it finds $\pthKt$
through \ex{6.18*}.

\noindent
{\bf Acknowledgements.}
I am very thankful to D.~Thurston and A.~Vaintrob for discussing this
work. I am especially indebted to A.~Vaintrob for numerous
discussions of the properties of the space $\cD$ and to D.~Thurston for substantive
discussions of the conjecture and for explaining the effects of
unwheeling procedure at the diagrammatic level.
This work was supported by NSF Grants
DMS-0196235 and DMS-0196131.

\appendix
\nsection{The 2-loop polynomial $\pthKt$ for knots with up to 7 crossing}

\begin{table}
\begin{center}
\begin{tabular}{||c|l|l||} \hline\hline
Knot & $\APKt$ & $12\pthKt$ \\
\hline\hline
$3_1$ & $t-1$ & $-t_1^2 t_2 + t_1^2$
\\ \hline
$4_1$ & $t^2 - 3t + 5$ & $0$
\\ \hline
$5_1$ & $t^2 - t + 1$ &
$2t_1^4 t_2^2 - 2t_1^4 t_2 + 2 t_1^4 - t_1^2 t_2 + t_1^2$
\\ \hline
$5_2$ & $2t - 3$ & $-13 t_1^2 t_2 + 9 t_1^2 + 6 t_1 - 12$
\\ \hline
$6_1$ & $-2t+5$ & $3t_1^2 t_2 - t_1^2 - 6t_1 + 24$
\\ \hline
$6_2$ & $-t^2 + 3t - 3$ & $-3t_1^4 t_2^2 +
3t_1^4 t_2 - t_1^4 - 6 t_1^3 - 11 t_1^2 t_2 + 15t_1^2$
\\ \hline
$6_3$ & $t^2 - 3t + 5$ & $0$
\\ \hline
$7_1$ & $ t^3 - t^2 + t - 1$ &
$-3t_1^6 t_2^3 + 3t_1^6 t_2^2 - 3t_1^6 t_2 + 2t_1^4 t_2^2 + 3t_1^6
 - 2 t_1^4 t_2 + 2 t_1^4 - t_1^2 t_2 + t_1^2$
\\ \hline
$7_2$ & $3t - 5$ & $-58t_1^2 t_2+ 36 t_1^2 + 36 t_1 - 84$
\\ \hline
$7_3$ & $2 t^2 - 3 t + 3$ & $-25 t_1^4 t_2^2 + 25 t_1^4 t_2 - 17 t_1^4 + 7
t_1^2 t_2 - 12 t_1^3 + t_1^2 - 6 t_1 + 12$
\\ \hline
$7_4$ & $4t-7$ & $ 136 t_1^2 t_2 - 80 t_1^2 - 96 t_1 + 240$
\\ \hline
$7_5$ & $2t^2 - 4t + 5$ &
$41 t_1^4 t_2^2 - 33 t_1^4 t_2 - 16 t_1^3 t_2 + 17 t_1^4 + 12 t_1^2 t_2
+32 t_1^3 + 4 t_1^2 - 14 t_1 + 36$
\\ \hline
$7_6$ & $-t^2 + 5t - 7$ &
$-7 t_1^4 t_2^2 + 5 t_1^4 t_2 + 10 t_1^3 t_2 - t_1^4 - 20 t_1^3
- 98 t_1^2 t_2  + 64 t_1^2 + 50 t_1 - 108$
\\ \hline
$7_7$ & $t^2 - 5t + 9$ & $ - 5 t_1^2 t_2 + t_1^2 + 12 t_1 -48 $
\\ \hline\hline
\end{tabular}
\end{center}
\caption{The Alexander polynomial $\APKt$ and the
2-loop polynomial $12\pthKt$ presented by monomials in funamental
domains}
\label{t6.1}
\end{table}

\begin{table}
\begin{center}
\begin{tabular}{||c|l|l||}
\hline\hline
Knot & $\APKt$ & $12\pthKt$
\\ \hline\hline
$3_1$ & $u-1$ &
$u_1^2 - 3u_2 u_2 - 2u_1 - 6$
\\ \hline
$4_1$ & $u^2 - 3u + 3$ & $0$
\\ \hline
$5_1$ & $u^2 - u - 1$ &
$2u_1^4 - 10 u_1^3 u_2 + 10 u_1^2 u_2^2 - 4 u_13 + 10 u_1^2 u_2 - 23 u_1^2
+ 53 u_1 u_2 $
\\
& & $+ 26 u_1 +66$
\\ \hline
$5_2$ & $2u-3$ &
$ 9u_1^2 - 31u_1 u_2 - 12 u_1 - 66$
\\ \hline
$6_1$ & $-2u+5$ &
$-u_1^2 + 5u_1 u_2 - 4 u_1 + 30$
\\ \hline
$6_3$ & $u^2 - 3u + 3$ & $0$
\\ \hline
$7_1$ & $u^3 - u^2 - 2u + 1$ &
$3u_1^6 - 21 u_1^5 u_2 + 42 u_1^4 u_2^2 - 21 u_1^3 u_2^3
- 6 u_1^5 + 21 u_1^4 u_2 - 52 u_1^4$
\\
&& $ + 215 u_1^3 u_2 - 152 u_1^2 u_2^2 +
62 u_1^3 - 26 u_1^2 u_2 + 268 u_1^2 - 358 u_1 u_2$
\\
&& $ - 64 u_1 - 276$
\\ \hline
$7_2$ & $3u-5$ & $36u_1^2 - 130u_1 u_2 - 36 u_1 - 300$
\\ \hline
$7_3$ & $2u^2 - 3u - 1$ &
$-17u_1^4 + 93u_1^3 u_2 - 109 u_1^2 u_2^2 + 38u_1^3 - 121 u_1^2 u_2
+ 221 u_1^2$
\\
&& $ - 559 u_1 u_2 - 314 u_1 - 702$
\\ \hline
$7_4$ & $4u-7$ & $-80u_1^2 + 296 u_1 u_2 + 64 u_1 + 720$
\\ \hline
$7_5$ & $2u^2 - 4u + 1$ &
$17 u_1^4 - 101 u_1^3 u_2 + 141 u_1^2 u_2^2 - 50 u_1^3 + 165 u_1^2 u_2
- 200 u_1^2$
\\
&& $ + 624 u_1 u_2 + 392 u_1 + 672$
\\ \hline
$7_6$ & $-u^2 + 5u - 5$ &
$-u_1^4 + 9 u_1^3 u_2 - 19 u_1^2 u_2^2 - 6 u_1^3 + 17 u_1^2 u_2 + 64 u_1^2
- 194 u_1 u_2$
\\
&& $ - 16 u_1 - 324$
\\ \hline
$7_7$ & $u^2 - 5u + 7$ &
$u_1^2 - 7 u_1 u_2 + 10 u_1 - 54$
\\ \hline\hline
\end{tabular}
\end{center}
\caption{The Alexander polynomial $\APKt$ and the
2-loop polynomial $12\pthKt$ expressed in terms of symmetric polynomials $u$
and $u_1$, $u_2$}
\label{t6.2}
\end{table}

Here are the results of calculating the polynomials $\pthKt$ for the first few
knots (with up to 7 crossings). We present these results in two different ways.
First, as we know, $\pth(\cK)\in\Qalghzo$ and relations\rx{6.8} come from the
symmetry $\Sgri{D_1}$. More explicitly, $\Hoz{D_1}$ looks like $su(3)$ root
lattice with elements $f_1$ and $f_2$ corresponding to the simple roots. The
symmetry group $\Sgri{D_1}$ is the symmetry of this lattice (which preserves
the origin). So instead of writing the whole polynomial $\pthKt$ we may list
just the monomials belonging to a fundamental domain of $\Sgri{D_1}$. From our
$su(3)$ lattice description it is easy to see that we may choose a fundamental
domain to include the monomials
\qq
t_1^{m_1} t_2^{m_2},\qquad m_1,m_2\geq 0,\; m_1 \geq 2 m_2.
\label{6.20*}
\qqq
Then the other monomials will be determined by the symmetries\rx{6.8}.
Similarly, in view of the symmetry\rx{4.5*} it is enough to list only the
monomials of $\APKt$ with non-negative powers of $t$. Thus in \tb{t6.1} we
present the `fundamental domain' parts of the Alexander polynomial $\APKt$ and
(scaled) 2-loop polynomial $12\pthKt$.

An alternative way of describing $\pthKt$ comes from the observation that the ring
of Laurent polynomials with the symmetries\rx{6.8} can be written as
$\IQ[u_1,u_2]$, where
\qq
u_1(t_1,t_2) &=&
t_1 + t_1^{-1} + t_2 + t_2^{-1} + t_1 t_2 + t_1^{-1} t_2^{-1},
\nonumber\\
u_2(t_1,t_2) & = &
t_1^2 t_2 + t_1^{-2} t_2^{-1} + t_1 t_2^{2} + t_1^{-1} t_2^{-2}
+ t_1 t_2^{-1} + t_1^{-1} t_2.
\label{6.21*}
\qqq
So in \tb{t6.2} we present the expressions for the Alexander polynomial $\APKt$
in terms of $u=t+t^{-1}$ and for the (scaled) 2-loop polynomial $12\pthKt$ in
terms of $u_1$ and $u_2$.

\begin{remark}
\rm
If $\cK\p$ is the mirror image of $\cK$, then
$\pth(\cK\p;\tot) = -\pthKt$, hence $\pthKt=0$ for amphicheiral knots.
\end{remark}
\begin{remark}
\rm
As we see, experimental evidence suggests that
\qq
12\pthKt\in\ZZtot.
\label{6.19*}
\qqq
\end{remark}
\begin{remark}
\rm
The degree of the Alexander polynomial is bounded by the genus of the knot
$g(\cK)$
\qq
\deg \APKt \leq g(\cK),
\label{6.22*}
\qqq
In view of the symmetries\rx{6.8} (which come from $\Sgri{D_1}$), the reasonable
measure of the degree of $\pthKt$ is the $t_1$ degree of its fundamental domain
part. Let us denote it simply as $\deg\pthKt$. Then \tb{t6.1} suggests a
similar inequality
\qq
\deg\pthKt\leq 2g(\cK),
\label{6.23*}
\qqq

\end{remark}


\end{document}